\documentclass[twoside,11pt]{article}
\title{} \author{} \date{}
\usepackage{latexsym,amssymb,times}%,showkeys}
\newtheorem{te}{Theorem}[section]

\newtheorem{cor}[te]{Corollary}
\newtheorem{fac}[te]{Fact}
\newtheorem{lem}[te]{Lemma}

\def\dok{\noindent{\bf Proof. }}
\def\kdok{\hfill $\Box$ \par \vspace*{2mm} }
\def\a{\alpha}

\def\o{\omega}

\def\f{\varphi}
\def\k{\kappa}
\def\p{\psi}
\def\r{\rho}

\def\t{\tau}

\def\L{\Lambda}

\def\B{\mathbb B}
\def\N{\mathbb N}
\def\P{\mathbb P}
\def\X{\mathbb X}
\def\Y{\mathbb Y}
\def\S{\mathbb S}
\def\Q{\mathbb Q}
\def\A{{\mathcal A}}
\def\D{{\mathcal D}}
\def\I{{\mathcal I}}
\def\CS{{\mathcal S}}
\def\CL{{\mathcal L}}

\def\la{\langle}
\def\ra{\rangle}
\def\dom{\mathop{\rm dom}\nolimits}
\def\ran{\mathop{\rm ran}\nolimits}
\def\Lev{\mathop{\rm Lev}\nolimits}
\def\Sl{\mathop{\rm Sl}\nolimits}
\def\Fn{\mathop{\rm Fn}\nolimits}
\def\Ip{\mathop{\rm Ip}\nolimits}
\def\height{\mathop{\rm ht}\nolimits}
\def\Emb{\mathop{\rm Emb}\nolimits}
\def\sq{\mathop{\rm sq}\nolimits}
\def\ro{\mathop{\rm ro}\nolimits}
\def\Fin{\mathop{\rm Fin}\nolimits}
\def\Scatt{\mathop{\rm Scatt}\nolimits}
\begin{document}
\thispagestyle{plain}
\begin{center}
           {\large \bf \uppercase{Copies of the Random Graph}}
\end{center}
\begin{center}
{\bf Milo\v s S.\ Kurili\'c\footnote{Department of Mathematics and Informatics, Faculty of Science, University of Novi Sad,
              Trg Dositeja Obradovi\'ca 4, 21000 Novi Sad, Serbia.
              email: milos@dmi.uns.ac.rs}
and Stevo Todor\v cevi\'c\footnote{Institut de Math\'ematique de Jussieu
              (UMR 7586) Case 247, 4 Place Jussieu, 75252 Paris Cedex, France
              and
              Department of Mathematics, University of Toronto,
              Toronto, Canada M5S 2E4.
              email: stevo@math.univ-paris-diderot.fr and stevo@math.toronto.edu}}
\end{center}
\begin{abstract}
\noindent
Let $\la R, \sim \ra$ be the Rado graph, $\Emb (R)$ the monoid of its self-embeddings, $\P (R)=\{ f(R): f\in \Emb (R)\}$ the set of copies
of $R$ contained in $R$, and ${\mathcal I}_R$ the ideal of subsets of $R$ which do not contain a copy of $R$.
We consider the poset $\la \P (R ), \subset\ra$, the algebra $P (R)/{\mathcal I _R}$, and the inverse of the right Green's pre-order on $\Emb (R)$,
and show that these pre-orders are forcing equivalent to a two step iteration of the form $\P \ast \pi$,
where the poset $\P$ is similar to the Sacks perfect set forcing: adds a generic real, has the $\aleph _0$-covering property and, hence, preserves $\o _1$, has the Sacks property and does not produce splitting reals, while $\pi$ codes an $\o$-distributive forcing.
Consequently, the Boolean completions of these four posets are isomorphic and the same holds for each countable graph containing a copy of the Rado graph.
\vspace{1mm}
\end{abstract}
\section{Introduction}\label{S1}
In this paper we continue the investigation of
%posets of isomorphic substructures of ultrahomogeneous relational structures. Generally speaking, we investigate
the partial orderings of the form $\la \P (\X ), \subset \ra$, where $\X $ is an ultrahomogeneous relational structure and
$\P (\X )$ the set of domains of  substructures of $\X$ isomorphic to $\X$.
In particular, if $\X = \la X, \r \ra$ is a binary structure (that is $\r \subset X\times X$), then
$\P (\X )  =  \{ A\subset X : \la A, \r _A \ra \cong \la X, \r \ra\} $, where $\r _A = \r \cap (A\times A)$.
In the sequel, in order to simplify notation, instead of $\la \P (\X ), \subset \ra$ we will write $\P (\X )$ whenever the context admits.

This investigation is related to a coarse classification of relational structures.
Namely, the conditions $\P (\X )=\P (\Y )$,  $\P (\X ) \cong \P (\Y )$,  $\sq \P (\X )\cong \sq \P (\Y )$ and $\ro\sq \P (\X )\cong \ro \sq \P (\Y )$
(where $\sq \P$ denotes the separative quotient of a partial order $\P$ and $\ro \sq \P$ its Boolean completion) define different equivalence relations
(``similarities") on the class of relational structures and their interplay with the similarities defined by the conditions
$\X =\Y$, $\X\cong \Y$ and $\X\rightleftarrows \Y$ (equimorphism) was considered in \cite{Kdif}. It turns out that the similarity
defined by the condition $\ro\sq \P (\X )\cong \ro \sq \P (\Y )$ is implied by all the similarities listed above and, thus, provides
the coarsest among the mentioned classifications of relational structures. Since the posets of copies are always homogeneous,
the condition $\ro\sq \P (\X )\cong \ro \sq \P (\Y )$ is equivalent to the forcing equivalence of the posets $\P (\X )$
and $\P (\Y )$ (we will write $\P (\X )\equiv \P (\Y )$) and, for convenience, we will exploit this fact using the tools of set-theoretic forcing in our proofs.

This paper can also be regarded as a part of the investigation of the quotient algebras of the form
$P(\omega )/{\mathcal I}$, where ${\mathcal I}$ is an  ideal on $\o$.
Namely, by \cite{Ktow}, if $\X $ is a countable indivisible structure with domain $\o$, then the collection ${\mathcal I}_\X $ of subsets of $\o$ which do not contain
a copy of $\X$ is either the ideal of finite sets or a co-analytic tall ideal and the poset $\sq \P (\X )$ is isomorphic to a dense subset
of $(P(\omega )/{\mathcal I}_\X)^+$, which implies $\ro\sq \P (\X )\cong \ro (P(\omega )/{\mathcal I}_\X)^+ $.
So, since the structure considered in this paper, the Rado graph, $\langle R, \sim \rangle $, is indivisible, our results can be regarded as statements concerning the forcing related properties of the corresponding quotient algebra. Namely, if we call a graph scattered if it does
not contain a copy of the Rado graph, and if ${\mathcal I}_R$ denotes the ideal of scattered subgraphs of $R$, then
$$
\ro \sq \P (R)= \ro((P (R)/{\mathcal I _R})^+ ).
$$
As a consequence of the main result of  \cite{KurTod} we have
the following statement describing the forcing related properties of the poset of copies of the rational line, $\Q$, and the corresponding quotient
$P(\Q )/\Scatt$, where $\Scatt$ denotes the ideal of scattered suborders of $\Q$. Namely, if $\S$ denotes the Sacks perfect set forcing and sh$(\S )$
the size of the continuum in the Sacks extension, then we have
\begin{te}  \label{T4120}
For each countable non-scattered linear order $L$ and, in particular, for the rational line, the poset $\P (L)$ is forcing equivalent to
the two-step iteration
$$
\S \ast \pi,
$$
where $1_\S \Vdash `` \pi $ is a $\sigma$-closed forcing".
If the equality sh$(\S )=\aleph _1$ (implied by CH)
or PFA holds in the ground model, then the second iterand is forcing equivalent to the poset $(P(\o )/\Fin )^+$ of the
Sacks extension.
Consequently,
$$
\ro \sq \P (\Q ) \cong \ro ((P(\Q )/\Scatt )^+)\cong \ro (\S \ast \pi ).
$$
\end{te}
(We note that by \cite{Kur1} the poset of copies of a countable {\it scattered} linear order $L$
is forcing equivalent to a separative atomless $\omega _1 $-closed poset; thus, under CH, to $(P(\o )/\Fin )^+$ and then
$\ro\sq \P (L )\cong \ro (P(\o )/\Fin )^+$. The posets of copies of countable ordinals are described in \cite{Kurord}.)

In this paper we prove a similar statement for non-scattered graphs (that is, the graphs containing a copy of the Rado graph):
\begin{te}  \label{T4}
For each countable non-scattered graph $\la G , \sim \ra$ and, in particular, for the Rado graph, the poset $\P (G)$ is forcing equivalent to
the two-step iteration
$$\P \ast \pi, $$
where $1_\P \Vdash ``\pi \mbox{ is an }\o\mbox{-distributive forcing}"$ and the poset $\P$ is similar to the Sacks forcing: adds a generic real,
has the $\aleph _0$-covering property (thus preserves $\o _1$), has the Sacks property and does not produce splitting reals.
In addition,
$$
\ro \sq \P (G) \cong \ro(P (R)/{\mathcal I _R})^+ \cong \ro (\P \ast \pi )
$$
and these complete Boolean algebras are weakly distributive\footnote{A complete Boolean algebra $\B$ is called {\it weakly distributive}
(or {\it $(\o , \cdot , <\!\!\o)$-distributive}) iff for each cardinal $\k$ and each matrix
$[b_{n \a }: \langle n , \a \rangle \in \o \times \k ]$ of elements of ${\mathbb B}$ we have
$$
\textstyle
\bigwedge _{n \in \o }\;\;
\bigvee _{\a \in \k }\;\;
b_{n \a}
=
\bigvee _{s: \o \rightarrow [\k ]^{<\o }}\;\;
\bigwedge _{n \in \o }\;\;
\bigvee _{\a \in s(n ) }
b_{n \a} .
$$
}.
\end{te}
In fact, if $\la G , \sim \ra$ is a countable graph containing a copy of the Rado graph, then these two structures are equimorphic and, by \cite{Kdif}, forcing equivalent. So it is sufficient to prove the previous theorem assuming that $\la G , \sim \ra$ is the Rado graph.

Finally we note that the results of this paper are related to the investigation of the monoids of self- embeddings.
We recall that the right Green's pre-order $\preceq ^R$ on a monoid $\la M , \cdot , 1 \ra$ is defined by $x\preceq ^R y$ iff $x\cdot z =y$, for some $z$.
It is easy to check (see \cite{Kmon}) that the poset of copies $\P (\X )$ of a structure $\X$ is isomorphic to the antisymmetric quotient of the
pre-order $\la \Emb (\X ), (\preceq ^R )^{-1}\ra$ and, consequently, these pre-orders are forcing equivalent. Thus, by Theorem \ref{T4}, for the Rado graph we have
$\la \Emb (R ), (\preceq ^R )^{-1}\ra \equiv (\P \ast \pi )$ and the Boolean completion of the pre-order
$\la \Emb (R ), (\preceq ^R )^{-1}\ra $ is a weakly distributive complete Boolean algebra.
\section{Preliminaries}\label{S2}
First we introduce a convenient notation. If $\la G, \sim \ra$ is a graph (namely, if $\sim$ is a symmetric and irreflexive binary relation on the set $G$) and $K\subset H \in [G]^{<\o }$, let
$$
G^H_K \!:=\Big\{ v\in G \setminus H : \forall k\in K \, (v\sim k) \; \land \; \forall h\in H\setminus K \, (v\not\sim h)\Big\}.
$$
(Clearly, $G^\emptyset _\emptyset =G$.)

The object of our study is the Rado graph (the Erd\H{o}s-R\'enyi graph, the countable random graph)
introduced independently by Erd\H{o}s and R\'enyi \cite{Erdos2} and Rado \cite{Rado}. It is characterized
as the unique (up to isomorphism) countable graph $\langle R, \sim \rangle $ such that
\begin{equation}\label{EQ2740}
R^H_K\neq \emptyset, \mbox{ whenever } K\subset H\in [R]^{<\o }.
\end{equation}
Equivalently, the Rado graph can be characterized as
the unique countable ultrahomogeneous universal graph (see \cite{Hodg})
or as the Fra\"{\i}ss\'e limit of the amalgamation class of all finite graphs (see \cite{Fra}).
In addition, by \cite{Erdos2}, if a graph with countably many vertices is chosen at random,
by picking edges independently with probability $\frac{1}{2}$, then, with probability 1, the
obtained graph will be isomorphic to the Rado graph.
The Rado graph and several related structures
(for example the automorphism group and the endomorphism monoid of $\langle R, \sim \rangle $, various topologies on $R$ etc.)
were extensively explored (see the survey article \cite{Camer}). The following fact contains the basic properties of the Rado graph which will be used in the paper.
\begin{fac}     \label{T2600}
Let $\langle R,\sim \rangle$ be a Rado graph and $\P (R)$ the set of its copies. Then

(a) If $F$ is a finite subset of $R$, then $R\setminus F \in \P (R)$;

(b) If $\{ X_1,  \dots , X_k\} $ is a partition of $R$, then $X_i\in \P (R)$, for some $i\leq k$ (the Rado graph is a strongly indivisible structure);

(c) If $H$ is a finite subset of $R$, then $\{ H \} \cup \{ R^H_K : K\subset H\}$ is a partition of $R$ and $R^H_K\in \P (R)$, for each $K\subset H$.
\end{fac}
Concerning the order theoretic properties of the poset $\langle \P (R) ,\subset \rangle$ we note that it is a
homogeneous, atomless and chain complete suborder of the order $\langle [R]^\omega ,\subset \rangle$ having a largest element, $R$.
In addition, by \cite{KuMa}, it contains maximal antichains of size ${\mathfrak c}$,
$\aleph _0$ and $n$, for each positive integer $n$, and in \cite{KuKu} the order types of maximal chains in this poset
are characterized as the order types of sets of the form $K \setminus \{ \min K \}$, where $K$ is a compact subset of the real line having the minimum non-isolated.

The sets $R^H_K$ (the orbits of $R$) will play an important role in our constructions.
\begin{lem}     \label{T2620}
Let $H_1$ and $H_2$ be finite subsets of $R$,  $K_1\subset H_1$ and $K_2\subset H_2$. Then

(a) $R^{H_1}_{K_1}\cap R^{H_2}_{K_2}\neq \emptyset $ if and only if $H_1 \cap K_2 = H_2 \cap K_1$;

(b) $R^{H_1}_{K_1}\cap R^{H_2}_{K_2}\neq \emptyset$ implies that $R^{H_1}_{K_1}\cap R^{H_2}_{K_2}=R^{H_1\cup H_2}_{K_1\cup K_2}$;

(c) $R^{H_1}_{K_1}= R^{H_2}_{K_2}$ if and only if $H_1 = H_2$ and $K_1= K_2$;

(d) $R^{H_1}_{K_1}\subset R^{H_2}_{K_2}$ if and only if $H_1 \supset H_2$,  $K_1\supset K_2$ and $H_2\cap K_1 =K_2$.
\end{lem}
\dok
We prove (a) and (b) simultaneously.

$(\Rightarrow )$ Assuming that $v\in R^{H_1}_{K_1}\cap R^{H_2}_{K_2}$ we first show that $H_1 \cap K_2 \subset H_2 \cap K_1$. If $r\in H_1 \cap K_2 $,
then $r\in H_2$ and, since $v\in  R^{H_2}_{K_2}$, we have $v\sim r$. Now $r \not\in K_1$ would imply $r\in H_1 \setminus K_1$ and,
since $v\in R^{H_1}_{K_1}$, we would have $v\not\sim r$, which is not true. So $r\in K_1$ and we are done. The reversed inclusion has a symmetric proof.

$(\Leftarrow )$ Let $H_1 \cap K_2 = H_2 \cap K_1$. Since $H_1 =(H_1 \setminus H_2) \dot{\cup} (H_1 \cap H_2)$ and
$K_1 =(K_1 \setminus H_2) \dot{\cup} (K_1 \cap H_2)$ and the second partition refines the first we have
\begin{equation}\label{EQ2679}
R^{H_1}_{K_1}= R^{H_1 \setminus H_2}_{K_1 \setminus H_2} \cap R^{H_1 \cap H_2}_{K_1 \cap H_2}.
\end{equation}
Similarly, $H_2 =(H_2 \setminus H_1) \dot{\cup} (H_2 \cap H_1)$ and $K_2 =(K_2 \setminus H_1) \dot{\cup} (K_2 \cap H_1)$, thus
\begin{equation}\label{EQ2680}
R^{H_2}_{K_2}= R^{H_2 \setminus H_1}_{K_2 \setminus H_1} \cap R^{H_2 \cap H_1}_{K_2 \cap H_1}.
\end{equation}
Now since $\{ H_1 \setminus H_2, H_1 \cap H_2 , H_2 \setminus H_1\}$ is a partition of the set $H_1 \cup H_2$ and, by the assumption,
$\{ K_1 \setminus H_2, K_1 \cap H_2 , K_2 \setminus H_1\}$ a partition of the set $K_1 \cup K_2$ refining the mentioned partition of $H_1 \cup H_2$, we have
\begin{equation}\label{EQ2617}
R^{H_1}_{K_1}\cap R^{H_2}_{K_2}
= R^{H_1 \setminus H_2}_{K_1 \setminus H_2} \cap R^{H_1 \cap H_2}_{K_1 \cap H_2}\cap  R^{H_2 \setminus H_1}_{K_2 \setminus H_1}
=R^{H_1\cup H_2}_{K_1\cup K_2}.
\end{equation}
and by Fact \ref{T2600}(c), $R^{H_1\cup H_2}_{K_1\cup K_2}\neq\emptyset$ (moreover, this set is a copy of $R$).

(c) Let $R^{H_1}_{K_1}= R^{H_2}_{K_2}$. Suppose that $H_2 \setminus H_1\neq \emptyset$ and let $v\in H_2 \setminus H_1$.
If $v\in K_2$, then, since $v\not\in H_1$, there is $w\in R^{H_1}_{K_1}\cap R^{\{v\}}_{\;\emptyset}$, thus $w\not\sim v$. But
$w\in R^{H_2}_{K_2}$ and, since  $v\in K_2$ we have $w\sim v$, so we have a contradiction. Otherwise, if $v\in H_2 \setminus K_2$,
there is $w\in R^{H_1}_{K_1}\cap R^{\{v\}}_{\{v\}}$, thus $w\sim v$. But
$w\in R^{H_2}_{K_2}$ and, since  $v\not\in K_2$ we have $w\not\sim v$ and we get a contradiction again. Thus $H_2 \subset H_1$ and,
similarly, $H_1 \subset H_2$. So $H_1 = H_2$, which by (a) implies $H_1\cap K_2= H_1 \cap K_1$, that is $K_2=K_1$.

(d) If $R^{H_1}_{K_1}\subset R^{H_2}_{K_2}$, then by (b) we have $R^{H_1}_{K_1}= R^{H_1}_{K_1}\cap R^{H_2}_{K_2}=R^{H_1\cup H_2}_{K_1\cup K_2}$
which, by (c), implies $H_1 = H_1\cup H_2$ and $K_1 = K_1\cup K_2$, that is $H_1 \supset H_2$ and  $K_1\supset K_2$.
Hence $H_1 \cap K_2= K_2$ so, by (a), $H_2\cap K_1 =K_2$

If $H_1 \supset H_2$,  $K_1\supset K_2$ and $H_2\cap K_1 =K_2$, then, since $K_2 =H_1 \cap K_2 $, by (a) and (b) we have
$R^{H_1}_{K_1}\cap R^{H_2}_{K_2}=R^{H_1\cup H_2}_{K_1\cup K_2}=R^{H_1}_{K_1}$, that is $R^{H_1}_{K_1}\subset R^{H_2}_{K_2}$.
\kdok
The reader will notice that, by (c) and (d) of the previous lemma, the mapping $F:\Fn (\o ,2)\rightarrow \P(R)$ given by $F(\f )= R^{\dom \f}_{\f ^{-1}[\{ 1\}]}$
is an embedding of the Cohen poset $\la \Fn (\o ,2), \supset \ra$ into the poset $\la \P (R), \subset \ra$. But $F$ is not a dense embedding
(we recall that $\P (R)$ contains antichains of size ${\mathfrak c}$) and this fact does not
imply that the poset $\P (R)$ is forcing equivalent to the Cohen forcing.
\section{Labeling of the vertices of the Rado graph}\label{S3}
Let $\la R , \sim \ra$ be the Rado graph.  A {\it labeling} of $L\in \P (R)$ is a pair $\CL =\la \Pi , q \ra$, where

(L1) $\Pi =\{ L_n : n\in \o\}$ is a partition of the set $L$,

(L2) $q: \bigcup _{n\in \o }\{ n \} \times P(\bigcup _{i<n} L_i) \rightarrow L$ is a bijection,

(L3) $L_n =\{ q (n,K) : K\subset \bigcup _{i<n}L_i \}$, for each $n\in \o$,

(L4) $q (n,K) \in  L ^{\bigcup _{i<n}L_i}_K$, for each $n\in \o$ and each $K\subset \bigcup _{i<n}L_i$.

\noindent
Then, clearly, $L_0=\{q(0 ,\emptyset ) \}$,  $|L_0|=1$ and the sets $L_n$ are finite. More precisely, by (L3) we have $|L_n |=m_n$,
where the  integers $m_n$, $n\in \o$, are defined by: $m_0=1$ and $m_n =2^{\sum _{i<n}m_i}$, for $n>0$.
Thus $\la |L_n|:n\in \o \ra =\la 1,2,8, 2^{11}, \dots \ra$.
\begin{lem} \label{T2639}
Each copy $L$ of $R$ has infinitely many labelings.
\end{lem}
\dok
Let $\prec _0$ be a well ordering on $L$ such that $\la L, \prec_0\ra\cong \la \o , <\ra$, where $<$ is the natural ordering on $\o$;
in fact w.l.o.g.\ we can assume that $\la L, \prec _0 \ra = \la \o , <\ra$. By recursion we define a sequence $\la L_n : n\in \o\ra$ such that for each
$m, n\in \o$ we have

(i)  $\;L_n$ is a finite subset of $L$,

(ii)  $L_m \cap L_n =\emptyset $, if $m\neq n$,

(iii) $L_n = \{ \min  L ^{\bigcup _{i<n}L_i}_K : K\subset \bigcup _{i<n}L_i \}$.

\noindent
First, since $L^\emptyset _\emptyset =L $, the sequence $\la L_0\ra$, where $L_0=\{ 0 \} $ satisfies (i), (ii) and (iii).

If $n>0$ and if $\la L_i : i<n\ra$ is a sequence satisfying (i) - (iii), then  $\bigcup _{i<n}L_i$ is a finite subset of $L$
and, by Fact \ref{T2600}(c), for $K\subset \bigcup _{i<n}L_i$ we have  $L ^{\bigcup _{i<n}L_i}_K \neq \emptyset$. Thus we define
$L_n$ by (iii) and, since $(\bigcup _{i<n}L_i)\cap L ^{\bigcup _{i<n}L_i}_K=\emptyset$, the extended sequence $\la L_i : i<n+1 \ra$ satisfies (i) - (iii). The recursion works.

Suppose that there is $m\in L\setminus \bigcup _{n\in \o }L_n$. For $n\in \o$, since $m\not\in \bigcup _{i<n}L_i$, by Fact \ref{T2600}(c)
there is $K_n\subset \bigcup _{i<n}L_i $ such that $m\in L ^{\bigcup _{i<n}L_i}_{K_n}$ and, since $m\not\in L_n$, by (iii)
we have $m> \min  L ^{\bigcup _{i<n}L_i}_{K_n}\in L _n$. Thus for each $n\in \o$ there is $q\in L_n$ such that $m>q$ and, by (ii),
$m$ is greater than infinitely many natural numbers, which is impossible. So, $\Pi :=\{ L_n : n\in \o\}$ is a partition of the set $L$.

Let the mapping $q: \bigcup _{n\in \o }\{ n \} \times P(\bigcup _{i<n} L_i) \rightarrow L$  be defied by
$q(n,K)=\min  L ^{\bigcup _{i<n}L_i}_K$. Since $L=\bigcup _{n\in \o }L_n$ the mapping $q$ is a surjection.
If $q(n,K)=q(n',K')$, then, since  $q(n,K)\in L_n$, by (ii) we have $n=n'$, which by Fact \ref{T2600}(c) implies $K=K'$.
Thus $q$ is a bijection, (L3) and (L4) follow from (iii) and $\CL =\la \Pi , q \ra$ is a labeling of $L$ determined by the well ordering $<$.

Clearly, if $n\in \o$ and $\prec _n$ is a well ordering on $L= \o$ such that $\la L, \prec_n \ra$ $\cong \la \o , <\ra$ and $\min _{\prec _n}L =n$
then repeating the previous construction using $\prec _n$ instead of $\prec _0$ we obtain a labeling $\CL _n$ of $L$, for which we have $L_0=\{ n\}$.
So the labelings $\CL _n$, $n\in \o$, are different.
\kdok
For convenience, instead of $q(n,K)$ we will write $q ^{\bigcup _{i<n}L_i}_K$ and a labeling
will be denoted by
$$\textstyle
\Big\la \{ L_n : n\in \o \}, \{ q ^{\bigcup _{i<n}L_i}_K : n\in \o \land K\subset \bigcup _{i<n}L_i\}\Big\ra.
$$
\section{Copies with orbits refining maximal antichains}\label{S4}
\noindent
The following construction of copies of $R$ will be frequently used in the paper.
We note that if $A\in \P( R)$  and $K\subset H\in [A]^{<\o }$, then $A$ with the induced graph structure is a Rado graph,
and, clearly, $\P (A)=P (A)\cap \P (R)$ and $A^H_K = A\cap R^H_K$.
\begin{lem} \label{T2633}
If $A\in \P (R)$ and if $L=\bigcup _{n\in \o} L_n$, where for each $n\in \o$ we have
$L_n =\{ q ^{\bigcup _{i<n}L_i}_K : K\subset \bigcup _{i<n}L_i \}$ and
$q ^{\bigcup _{i<n}L_i}_K \in  A ^{\bigcup _{i<n}L_i}_K$, for all $K\subset \bigcup _{i<n}L_i$, then

(a) $L\in \P (A)$;

(b) $\la \{ L_n : n\in \o \}, \{ q ^{\bigcup _{i<n}L_i}_K : n\in \o \land K\subset \bigcup _{i<n}L_i\}\ra$ is a labeling of $L$;

(c)  If $\f (u,v,w)$ is a formula of the language of set theory, $\t $ is a $\P (R)$-name and $S_n\in \P (A)$, for $n \in \o $,
where $L \subset S_{n+1}\subset S_n $, for each $n\in \o$,  and if
\begin{equation}\label{EQ2691}\textstyle
\forall n \in \o \;\; \forall K\subset \bigcup _{i<n}L_i\;\; \exists a^n_K \;\; (S_n) ^{\bigcup _{i<n}L_i}_K\Vdash \f (\t , \check{n}, \check{a^n_K}),
\end{equation}
then for each $n \in \o$ and each $K\subset \bigcup _{i<n}L_i$ we have $L^{\bigcup _{i<n}L_i}_K\Vdash \f (\t , \check{n}, \check{a^n_K})$.
\end{lem}
\dok
(a) If $K\subset H $ are finite subsets of $L$, then $H\subset \bigcup _{i<n}L_i$ for some $n\in \o$ and
$q ^{\bigcup _{i<n}L_i}_K \in L\cap A^{\bigcup _{i<n}L_i}_K \subset L\cap R^H_K= L^H_K$. Thus the graph $L$ satisfies (\ref{EQ2740}).

(b) By the assumption, the mapping $q: \bigcup _{n\in \o }\{ n \} \times P(\bigcup _{i<n} L_i) \rightarrow L$ defined by $q(n,K)=q^{\bigcup _{i<n}L_i}_K$
is a surjection. If $q^{\bigcup _{i<n_1}L_i}_{K_1} = q^{\bigcup _{i<n_2}L_i}_{K_2}=:u$, then, since
$u\in L_{n_1}\cap A^{\bigcup _{i<n_2}L_i}_{K_2}\subset L_{n_1} \cap R\setminus \bigcup _{i<n_2}L_i$ we have $n_1 \geq n_2$ and, similarly,
$n_2 \geq n_1$, which gives $n_1=n_2$. By Fact \ref{T2600}(c), $K_1\neq K_2$ implies
$A^{\bigcup _{i<n_1}L_i}_{K_1}\cap A^{\bigcup _{i<n_1}L_i}_{K_2}=\emptyset$, thus $K_1= K_2$, and, hence, $q$ is an injection. This implies that
$\{ L_n :n\in \o\}$ is a partition of $L$ and conditions (L3) and (L4) are obviously satisfied.

(c) Since $L \subset S_n$ we have $L^{\bigcup _{i<n}L_i}_K \subset (S_n) ^{\bigcup _{i<n}L_i}_K$ and we apply (\ref{EQ2691}).
\kdok
\noindent
Roughly speaking, in order to provide condition (\ref{EQ2691})
it is sufficient that for each maximal antichain $\A _n$ in $\P (R)$ such that each $A\in \A _n$ forces
$\f (\t , \check{n}, \check{a})$, for some $a$, there is $S_n\in \P (S_{n-1})$ containing $\bigcup _{i<n}L_i$ and such that each orbit
$(S_n) ^{\bigcup _{i<n}L_i}_K$ is contained in some $A\in \A _n$. This follows from the following theorem, the main statement of this section.
\begin{te}   \label{T2635}
For each maximal antichain $\A$ in the poset $\P (R)$ and each finite set $F_0\subset R$
there is $S\in \P (R)$ such that $F_0\subset S$ and
\begin{equation}\label{EQ2692}
\forall H \subset F_0 \;\; \exists A \in \A \;\; S^{F_0}_H \subset A .
\end{equation}
\end{te}
The proof of Theorem \ref{T2635}, given at the end of the section, is based on the following three lemmas.
\begin{lem}     \label{T2621}
For each $p\in R$ and each finite (possibly empty) sets $F\subset R ^{\{ p \}}_{\{ p \}}$ and $G\subset R ^{\{ p \}}_{\;\;\emptyset }$ there is an isomorphism
$f:R ^{\{ p \}}_{\{ p \}}\rightarrow R ^{\{ p \}}_{\;\;\emptyset }$  satisfying
\begin{equation}\label{EQ2615}
\forall u,v \in R ^{\{ p \}}_{\{ p \}} \;\;\Big( u\sim f(v) \;\Leftrightarrow \; f(u)\sim v \Big),
\end{equation}
\begin{equation}\label{EQ2666}
F\cap f^{-1}[G] =\emptyset.
\end{equation}
\end{lem}
\dok
Let $\P$ be the set of all partial functions $\f$ from $R ^{\{ p \}}_{\{ p \}}$ to $R ^{\{ p \}}_{\;\;\emptyset }$ such that
for each $u,v\in \dom \f$ we have

(i) $u\neq v \Rightarrow \f (u)\neq \f (v)$,

(ii) $u\sim v \Leftrightarrow \f (u)\sim \f (v)$,

(iii) $u\sim \f (v) \Leftrightarrow \f (u)\sim v$,

(iv) $u\in F \Rightarrow \f (u)\not\in G$.
\vspace{2mm}

\noindent
{\bf Claim 1.}
$D_u=\{ \f \in \P : u\in \dom \f \}$, $u\in R ^{\{ p \}}_{\{ p \}}$, are dense sets in $\la \P , \supset\ra$.

\vspace{2mm}

\noindent
{\it Proof of Claim 1.}
Let $u_1\in R ^{\{ p \}}_{\{ p \}}$ and $\p \in\P \setminus D_{u_1}$, that is, $u_1\not\in \dom \p$. Let
$$
H=\{ u\in \dom \p : u\sim u_1 \} \mbox{ and }
L=\{ v\in \dom \p : \p (v)\sim u_1 \}.
$$
Since the sets $\{ p \}$, $\dom \p$ and $\ran \p$ are pairwise disjoint we choose
\begin{equation}\label{EQ2635}
v_1 \in ( R ^{\{ p \}}_{\;\;\emptyset } \cap R ^{\dom \p }_L \cap R^{\ran \p }_{\p [H]})\setminus G
\end{equation}
and show that $\f =\p \cup \{ \la u_1 , v_1 \ra \} \in \P$. Since $u_1\not\in \dom \p$, $\f$ is a function and, by
(\ref{EQ2635}), we have $v_1\not\in \ran \p$ which implies (i).

Since $\p \in \P$ for a proof that $\f$ satisfies (ii) we show that
\begin{equation}\label{EQ2636}
\forall u\in \dom \p \;\;(u\sim u_1 \Leftrightarrow \p (u)\sim v_1 ).
\end{equation}
Let $u\in \dom \p$.
If $u\sim u_1$, then $u\in H$ thus $\p (u)\in \p [H]$ and, by (\ref{EQ2635}), $\p (u) \sim v_1$.
If $u\not\sim u_1$, then $u\in \dom \p \setminus H$ and, hence, $\p (u)\in \ran \p \setminus \p [H]$, which, together with (\ref{EQ2635}), implies
$\p (u)\not\sim v_1$. So (\ref{EQ2636}) is true.

For a proof that $\f$ satisfies (iii) we show that
\begin{equation}\label{EQ2637}
\forall v\in \dom \p \;\;(u_1\sim \p (v) \Leftrightarrow v_1\sim v ).
\end{equation}
Let $v\in \dom \p$.
If $\p (v)\sim u_1$, then $v\in L$ and, by (\ref{EQ2635}), $v_1\sim v$.
If $\p (v)\not\sim u_1$, then $v\in \dom \p \setminus L$ and, by (\ref{EQ2635}), $v_1\not\sim v$. Thus (\ref{EQ2637}) is true.

Since $\p \in \P$ and, by (\ref{EQ2635}), $\f (u_1)=v_1\not\in G$, $\f$ satisfies (iv).
Thus $\f \in \P$ and, clearly $\p \subset \f \in D_{u_1}$. Claim 1 is proved.
\kdok
\noindent
{\bf Claim 2.}
$\Delta _v=\{ \f \in \P : v\in \ran \f \}$, $v\in R ^{\{ p \}}_{\;\; \emptyset}$, are dense sets in $\la \P , \supset\ra$.

\vspace{2mm}

\noindent
{\it Proof of Claim 2.}
Let $v_1\in R ^{\{ p \}}_{\;\;\emptyset }$ and $\p \in\P \setminus \Delta _{v_1}$, that is, $v_1\not\in \ran \p$. Let
\begin{equation}\label{EQ2667}
H=\{ v\in \ran \p : v\sim v_1 \} \mbox{ and }
L=\{ u\in \dom \p : u\sim v_1 \}.
\end{equation}
Since the sets $\{ p \}$, $\dom \p$ and $\ran \p$ are pairwise disjoint we choose
\begin{equation}\label{EQ2638}
u_1 \in ( R ^{\{ p \}}_{ \{ p \} } \cap R ^{\dom \p }_{\p ^{-1}[H]} \cap R^{\ran \p }_{\p [L]})\setminus F
\end{equation}
and show that $\f =\p \cup \{ \la u_1 , v_1 \ra \} \in \P$. By (\ref{EQ2638}) we have $u_1\not\in \dom \p$ so $\f$ is a function and,
since $v_1\not\in \ran \p$, $\f $ satisfies (i).

Since $\p \in \P$, for a proof that $\f$ satisfies (ii) it remains to be shown that (\ref{EQ2636}) holds.
Let $u\in \dom \p$.
If $u\sim u_1$, then, by (\ref{EQ2638}), $u\in \p ^{-1}[H]$ thus $\p (u)\in H$ and, hence, $\p (u) \sim v_1$.
If $u\not\sim u_1$, then $u\in \dom \p \setminus \p ^{-1}[H]$ and, hence, $\p (u)\in \ran \p \setminus H$, which, by (\ref{EQ2667}), implies
$\p (u)\not\sim v_1$ and (\ref{EQ2636}) is true.

For a proof of (iii) we verify (\ref{EQ2637}).
Let $v\in \dom \p$.
If $v\sim v_1$, then $v\in L$ and, hence, $\p (v)\in \p [L]$ so, by (\ref{EQ2638}), $u_1\sim \p (v)$.
If $v\not\sim v_1$, then $v\in \dom \p \setminus L$ and, hence, $\p (v)\in \ran \p \setminus \p [L]$
and,  by (\ref{EQ2638}), $u_1\not\sim \p (v)$. Thus (\ref{EQ2637}) is true.

Since $\p \in \P$ and, by (\ref{EQ2638}), $u_1\not\in F$, $\f$ satisfies (iv).

Thus $\f \in \P$ and, clearly $\p \subset \f \in \Delta _{v_1}$. Claim 2 is proved.
\kdok

\noindent
By Claims 1, 2 and the Rasiowa-Sikorski theorem there is a filter ${\mathcal G}$ in the poset $\la \P ,\supset \ra$ intersecting the sets
$D_u$, $u\in R ^{\{ p \}}_{\{ p \}}$, and $\Delta _v$, $v\in R ^{\{ p \}}_{\;\; \emptyset}$. Thus
$f=\bigcup _{\f \in {\mathcal G}}\f \subset R ^{\{ p \}}_{\{ p \}} \times  R ^{\{ p \}}_{\;\; \emptyset}$ and
$\dom f = R ^{\{ p \}}_{\{ p \}}$ and $\ran f=R ^{\{ p \}}_{\;\; \emptyset}$. So, since ${\mathcal G}$ is a set of compatible functions,
$f$ is a surjection from  $R ^{\{ p \}}_{\{ p \}}$ onto $R ^{\{ p \}}_{\;\; \emptyset}$.
By (i) $f$ is an injection, by (ii) it is a graph-isomorphism, by (iii) satisfies (\ref{EQ2615}) and, by (iv), satisfies (\ref{EQ2666}).
\kdok

Let $p\in R$ and let $F\subset R ^{\{ p \}}_{\{ p \}}$ and $G\subset R ^{\{ p \}}_{\;\;\emptyset}$ be finite (possibly empty) sets.
A set $C\subset R ^{\{ p \}}_{\{ p \}}$ will be called {\it $(p,F,G)$-extendible} iff there is a set $C' \subset R ^{\{ p \}}_{\;\;\emptyset}$ such that
$F\cup G \subset C\cup \{ p \} \cup C'\in \P (R)$.
(Then, by Fact \ref{T2600}(c), $C$ and $C'$ are copies of $R$.) $(p,\emptyset ,\emptyset )$-extendible copies will be called $p$-extendible.

For $F=G=\emptyset$, the following statement shows that there is a copy  $B\subset R ^{\{ p \}}_{\{ p \}} $ such that the set of $p$-extendible copies
is dense below $B$. Moreover we have
\begin{lem}     \label{T2619}
Let $p\in R$, let $F\subset R ^{\{ p \}}_{\{ p \}}$ and $G\subset R ^{\{ p \}}_{\;\;\emptyset }$ be finite (possibly empty) sets and
$f:R ^{\{ p \}}_{\{ p \}}\rightarrow R ^{\{ p \}}_{\;\;\emptyset }$ an isomorphism satisfying (\ref{EQ2615}) and (\ref{EQ2666}).
Then there is a copy $B\in \P (R ^{\{ p \}}_{\{ p \}} )$ such that $F \cup f^{-1}[G] \subset B$ and that for each set $A$ satisfying
\begin{equation}\label{EQ2669}
F \cup f^{-1}[G] \subset A\in \P (B)
\end{equation}
there are sets  $A_0$ and $A_1$ such that
\begin{equation}\label{EQ2670}
A_0\cup A_1\subset A \;\; \land  \;\; A_0 \cap A_1 =\emptyset  \;\; \land \;\;  F\subset A_0  \;\; \land  \;\; G \subset f[A_1] ,
\end{equation}
\begin{equation}\label{EQ2633}
A_0 \cup \{ p \}\cup f[A_1]\in \P (R).
\end{equation}
Thus $A_0$ is a $(p,F,G)$-extendible copy contained in $ A$ .
\end{lem}
\dok
Let $B=\bigcup _{n\in \o}L_n $, where $L_0=F$, $L_1 =f^{-1}[G]$ and, for $n\geq 1$,
\begin{equation}\label{EQ2668}\textstyle
L_{n+1}=\{ q ^{\bigcup _{k\leq n}L_k}_K \!\!: K\subset \bigcup _{k\leq n}L_k \} , \mbox{ where}
\end{equation}
\begin{equation}\label{EQ2632}\textstyle
q ^{\bigcup _{k\leq n}L_k}_K \in  R ^{\{ p \}}_{\{ p \}} \cap R ^{\bigcup _{k\leq n}L_k}_K \cap R ^{f[\bigcup _{k\leq n}L_k]}_{f[K]},
\mbox{ for each }K\subset \bigcup _{k\leq n}L_k.
\end{equation}
Then $F \cup f^{-1}[G] \subset B\subset R ^{\{ p \}}_{\{ p \}}$ and, as in Lemma \ref{T2633}, we show that $B\in \P (R)$.

Let $A$ be a set satisfying (\ref{EQ2669}). We will construct sets $A_0$ and $A_1$ satisfying (\ref{EQ2670}) and (\ref{EQ2633}).

First by recursion we construct finite sets $S_{i,j}\subset \o \setminus \{ 0,1 \}$, for $2\leq i <\o$ and $j\in \{ 0,1 \}=2$, and $a_n\in A $, for
$n\in \bigcup _{2\leq i< \o}\bigcup _{j<2}S_{i,j}$, such that
\begin{itemize}
\item[(i)] $a_n\in A \cap L_n$, for $n\in \bigcup _{2\leq i< \o}\bigcup _{j<2}S_{i,j}$,

\item[(ii)] $\la i, j\ra <_{lex }\la i_1, j_1\ra$ implies $S_{i,j}<S_{i_1 , j_1}$ (that is, $\max S_{i,j}< \min S_{i_1 , j_1}$ ),

\item[(iii)] For each $i_0 \geq 2 $, each $K\subset L_0 \cup L_1 \cup A_{<i_0}$
(where, for simplicity, we define $S_{<i_0}:= \bigcup _{2\leq i<i_0 }\bigcup _{j<2}S_{i,j}$ and $A_{<i_0}:= \{ a_n : n\in S_{<i_0}\}$,
thus $S_{<2}$ $=A_{<2}=\emptyset$)  and each $j_0 <2$
there is $n$ such that
\begin{equation}\label{EQ2671}
n\in S_{i_0,j_0},
\end{equation}
\begin{equation}\label{EQ2672}
a_n \in R^{L_0 \cup L_1 \cup A_{<i_0} }_K .
\end{equation}
\end{itemize}
{\bf Claim 0.} The recursion works.

\vspace{2mm}
\noindent
{\it Proof of Claim 0.}
Let $i_0\geq 2$ and let $\la S_{i,j} :2\leq i<i_0 \land j<2 \ra$ and $\la a_n: n \in S_{<i_0} \ra$
satisfy conditions (i) - (iii).
Let $k_{i_0}=|L_0 \cup L_1 \cup S_{<i_0} |$ and let us fix an enumeration
\begin{equation}\label{EQ2673}\textstyle
P (L_0 \cup L_1 \cup A_{<i_0} )=\{ K_r : r<2^{k_{i_0} }\} .
\end{equation}
First we define $S_{i_0 ,0}$. Let $m_{i_0 ,0} =\max( \{ 0,1 \} \cup S_{<i_0} )$.
Since $A\in \P (R)$, for each $r<2^{k_{i_0} }$ the set  $A^{L_0 \cup L_1 \cup A_{<i_0} }_{K_r}$
is infinite and, by (\ref{EQ2669}), intersects infinitely many sets $L_n$.
So, for $r<2^{k_{i_0}}$ let
$$
                 n^r_{i_0,0}= \left\{ \begin{array}{ll}
                                    \min \{ n>m_{i_0,0}\!: A^{L_0 \cup L_1 \cup A_{<i_0} }_{K_0} \cap L_n \neq \emptyset \} &
                                    \mbox{ if } r=0,\\
                                    \min \{ n>n^{r-1}_{i_0 ,0}: A^{L_0 \cup L_1 \cup A_{<i_0} }_{K_r} \cap L_n \neq \emptyset \} &             \mbox{ if }r>0 ,
                                     \end{array}
                             \right.
$$
let us define $S_{i_0 ,0}=\{ n^r_{i_0 ,0 }: r< 2^{k_{i_0}}\}$ and for $r< 2^{k_{i_0}}$ let us choose
\begin{equation}\label{EQ2674}\textstyle
a_{n^r_{i_0,0}}\in A^{L_0 \cup L_1 \cup A_{<i_0} }_{K_r} \cap L_{n^r_{i_0,0}}.
\end{equation}
Now we define $S_{i_0 ,1}$. Let $m_{i_0 ,1} =\max( \{ 0,1 \} \cup S_{<i_0} \cup S_{i_0 ,0})$.
For $r<2^{k_{i_0}}$ let
$$
                 n^r_{i_0,1}= \left\{ \begin{array}{ll}
                                    \min \{ n>m_{i_0,1}\!: A^{L_0 \cup L_1 \cup A_{<i_0} }_{K_0} \cap L_n \neq \emptyset \} &
                                    \mbox{ if } r=0,\\
                                    \min \{ n>n^{r-1} _{i_0 ,1}: A^{L_0 \cup L_1 \cup A_{<i_0} }_{K_r} \cap L_n \neq \emptyset \} &             \mbox{ if }r>0 ,
                                     \end{array}
                             \right.
$$
let us define $S_{i_0 ,1}=\{ n^r_{i_0 ,1 }: r< 2^{k_{i_0}}\}$ and for $r< 2^{k_{i_0}}$ let us choose
\begin{equation}\label{EQ2675}\textstyle
a_{n^r_{i_0,1}}\in A^{L_0 \cup L_1 \cup A_{<i_0} }_{K_r} \cap L_{n^r_{i_0,1}}.
\end{equation}
By (\ref{EQ2673}), (\ref{EQ2674}) and (\ref{EQ2675}), the extended sequences
$\la S_{i,j} :2\leq i< i_0 +1 \land j<2 \ra$ and $\la a_n: n \in S_{<i_0 +1}\ra$
satisfy conditions (i) and (iii). By the construction we have $S_{< i_0}<S_{i_0 ,0}<S_{i_0 ,1}$ and (ii) is true as well. The recursion works indeed.
\kdok
\noindent
Now we define the sets $A_0$ and $A_1$ by:
\begin{equation}\label{EQ2676}%\textstyle
A_0= L_0 \cup \{ a_n : n\in \bigcup _{2\leq i< \o }S_{i,0}\}  \;\mbox{ and }\; A_1=L_1 \cup \{ a_n : n\in \bigcup _{2\leq i< \o }S_{i,1}\}.
\end{equation}
By (\ref{EQ2669}) and (i) we have $A_0 \cup A_1 \subset A$. By (\ref{EQ2666}) we have $L_0 \cap L_1 = F \cap f^{-1}[G]=\emptyset $, which,
together with (i), (ii) and (\ref{EQ2676}), implies $A_0 \cap A_1=\emptyset$.
By (\ref{EQ2676}) we have $F= L_0  \subset A_0$ and $f^{-1}[G] = L_1 \subset A_1$ so $G\subset f[A_1]$ and (\ref{EQ2670}) is true.

We prove (\ref{EQ2633}) showing that the set $A_0 \cup \{ p \}\cup f[A_1]$ satisfies (\ref{EQ2740}). Let
\begin{equation}\label{EQ2677}\textstyle
K_0 \subset H_0\in [A_0 ]^{<\o } \;\;\mbox{ and }\;\; K_1 \subset H_1 \in [f[A_1]]^{<\o }.
\end{equation}
Since  $f^{-1}[K_1]\subset f^{-1}[H_1]\subset A_1$, by (\ref{EQ2676}) there is $i_0> 2$ such that
\begin{equation}\label{EQ2621}\textstyle
K _0 \cup f^{-1}[K _1]\subset H_0 \cup f^{-1}[H_1]\subset L_0 \cup L_1 \cup A_{<i_0} .
\end{equation}
{\bf  Claim 1.}
$A_0 \cap R^{\{ p \} \cup H_0 \cup H_1}_{\{ p \} \cup K_0 \cup K_1} \neq \emptyset $.

\vspace{2mm}

\noindent
{\it  Proof of Claim 1.}
By (\ref{EQ2621}) and (iii) there is $n\in S_{i_0,0}$ such that
\begin{equation}\label{EQ2619}
a_n \in R^{L_0 \cup L_1 \cup A_{<i_0} }_{K_0 \cup f^{-1}[K_1]}.
\end{equation}
{\bf  Subclaim 1.1} $a_n \in R^{H_0}_{K_0}$.

\vspace{2mm}
\noindent
{\it  Proof of Subclaim 1.1}
For $u\in K_0$, by (\ref{EQ2619}) we have $a_n \sim u$.
For $u\in H_0 \setminus K_0$, by (\ref{EQ2621}) and (\ref{EQ2619}) and since $H_0 \cap f^{-1}[K_1]=\emptyset$, we have $a_n \not\sim u$.
\kdok
\noindent
{\bf  Subclaim 1.2} $a_n \in R^{H_1}_{K_1}$.

\vspace{2mm}
\noindent
{\it  Proof of Subclaim 1.2}
By the definition of $B$ and since $a_n\in A\setminus (L_0 \cup L_1) \subset B$ there are $n_0\geq 2$ and $K\subset \bigcup _{k\leq n_0}L_k$ such that
\begin{equation}\label{EQ2620}
a_n = q ^{\bigcup _{k\leq n_0}L_k}_K \in R ^{\{ p \}}_{\{ p \}} \cap R ^{\bigcup _{k\leq n_0}L_k}_K \cap R ^{f[\bigcup _{k\leq n_0}L_k]}_{f[K]}.
\end{equation}
Thus $a_n \in L_{n_0 +1}$ which by (i) implies $n= n_0 +1$ and, since $n\in S_{i_0 , 0}$, by (ii) we have
$S_{<i_0} <\{ n_0 +1\}$ and by (i),
\begin{equation}\label{EQ2624}\textstyle
L_0 \cup L_1 \cup A_{<i_0} \subset \bigcup _{k\leq n_0}L_k.
\end{equation}
By (\ref{EQ2619}), (\ref{EQ2620}) and Lemma \ref{T2620}(a) we have
$$
(L_0 \cup L_1 \cup A_{<i_0} ) \cap (\{ p \} \cup K \cup f[K])
= (\{ p \} \cup  \bigcup _{k\leq n_0}L_k \cup f[\bigcup _{k\leq n_0}L_k])\cap (K_0 \cup f^{-1}[K_1])
$$
which, by  (\ref{EQ2621}) and (\ref{EQ2624}), gives
\begin{equation}\label{EQ2625}\textstyle
(L_0 \cup L_1 \cup A_{<i_0} ) \cap  K = K_0 \cup f^{-1}[K_1].
\end{equation}
For $u\in K_1$ we have $f^{-1}(u)\in f^{-1}[K_1]$ and, by (\ref{EQ2625}),
$f^{-1}(u)\in K$ that is $u\in f[K]$ and, by (\ref{EQ2620}) we have $a_n \sim u$.

For $u\in H_1 \setminus K_1$ by (\ref{EQ2677}) we have $f^{-1}(u)\in f^{-1}[H_1]\setminus f^{-1}[K_1]\subset A_1$, which implies
$f^{-1}(u)\not\in K_0$. Thus
\begin{equation}\label{EQ2626}
f^{-1}(u)\not\in K_0 \cup f^{-1}[K_1].
\end{equation}
Since $f^{-1}(u)\in f[H_1]$, by (\ref{EQ2621})
$f^{-1}(u) \in L_0 \cup L_1 \cup  A_{<i_0} $ so, by (\ref{EQ2626}) and (\ref{EQ2625}),
$f^{-1}(u)\not\in K$ and, hence
\begin{equation}\label{EQ2627}
u\not\in f[K].
\end{equation}
By (\ref{EQ2621}) and (\ref{EQ2624}) we have $f^{-1}(u) \in\bigcup _{k\leq n_0}L_k$ which implies $u \in f[\bigcup _{k\leq n_0}L_k]$ so, by
(\ref{EQ2627}) and (\ref{EQ2620}), $a_n \not\sim u$.
\kdok
\noindent
Now, since $n\in S_{i_0,0}$ we have $a_n\in A_0$ and, by  (\ref{EQ2620}) and Subclaims 1.1 and 1.2,
$a_n\in A_0 \cap R^{\{ p \} \cup H_0 \cup H_1}_{\{ p \} \cup K_0 \cup K_1}$. Claim 1 is proved.
\kdok
\noindent
{\bf  Claim 2.}
$f[A_1] \cap R^{\{ p \} \cup H_0 \cup H_1}_{\;\;\emptyset\;\, \cup K_0 \cup K_1} \neq \emptyset $.

\vspace{2mm}
\noindent
{\it  Proof of Claim 2.}
By (\ref{EQ2621}) and (iii) there is $n\in S_{i_0,1}$ such that
\begin{equation}\label{EQ2628}
a_n \in R^{L_0 \cup L_1 \cup A_{<i_0} }_{K_0 \cup f^{-1}[K_1]}.
\end{equation}
{\bf  Subclaim 2.1} $f(a_n) \in R^{H_1}_{K_1}$.

\vspace{2mm}
\noindent
{\it  Proof of Subclaim 2.1}
For $u\in K_1$ we have $f^{-1}(u)\in f^{-1}[K_1]$ and, by (\ref{EQ2628}), $a_n \sim f^{-1}(u)$. Thus, since $f$ is an isomorphism,
$f(a_n)\sim u $.

For $u\in H_1 \setminus K_1$ we have  $f^{-1}(u)\in f^{-1}[H_1]\setminus f^{-1}[K_1]$
and, by (\ref{EQ2621}) and (\ref{EQ2628}) we have $a_n \not\sim f^{-1}(u)$. So, since $f$ is an isomorphism
$f(a_n)\not\sim u $.
\kdok
\noindent
{\bf  Subclaim 2.2} $f(a_n) \in R^{H_0}_{K_0}$.

\vspace{2mm}
\noindent
{\it  Proof of Subclaim 2.2}
By the definition of $B$ and since $a_n\in A_1 \setminus (L_0 \cup L_1)\subset B$ there are $n_0\geq 2$ and $K\subset \bigcup _{k\leq n_0}L_k$ such that
\begin{equation}\label{EQ2629}
a_n = q ^{\bigcup _{k\leq n_0}L_k}_K \in R ^{\{ p \}}_{\{ p \}} \cap R ^{\bigcup _{k\leq n_0}L_k}_K \cap R ^{f[\bigcup _{k\leq n_0}L_k]}_{f[K]}.
\end{equation}
Thus $a_n \in L_{n_0 +1}$ which by (i) implies $n= n_0 +1$ and, since $n\in S_{i_0 ,1}$, by (ii) we have
$S_{<i_0} <\{ n_0 +1\}$ and, by (i),
\begin{equation}\label{EQ2630}\textstyle
L_0 \cup L_1 \cup A_{<i_0} \subset \bigcup _{k\leq n_0}L_k.
\end{equation}
By (\ref{EQ2628}), (\ref{EQ2629}) and Lemma \ref{T2620}(a),
$$
(L_0 \cup L_1 \cup A_{<i_0} ) \cap (\{ p \} \cup K \cup f[K])
=(\{ p \} \cup  \bigcup _{k\leq n_0}L_k \cup f[\bigcup _{k\leq n_0}L_k])\cap (K_0 \cup f^{-1}[K_1])
$$
so, by  (\ref{EQ2621}) and (\ref{EQ2630}),
$(L_0 \cup L_1 \cup A_{<i_0})  \cap  K = K_0 \cup f^{-1}[K_1]$. Thus
\begin{equation}\label{EQ2631}\textstyle
f[L_0 \cup L_1 \cup A_{<i_0} ] \cap f[K] =
f[K_0 ]\cup K_1.
\end{equation}
Now, for $u\in K_0$ we have $f(u)\in f[K_0]$ and, by (\ref{EQ2631}), $f(u)\in f[K]$ which, by (\ref{EQ2629}) gives $a_n \sim f(u)$ and, by (\ref{EQ2615}),
$f(a_n)\sim u $.

For $u\in H_0 \setminus K_0$ we have  $f(u)\in f[H_0]\setminus f[K_0]$.
Since $u\in A_0$ and $f^{-1}[H_1]\subset A_1$ we have $u\not\in f^{-1}[H_1]$ and, hence $f(u)\not\in H_1\supset K_1$.
Thus $f(u)\not\in f[K_0]\cup K_1$.
Since $u\in H_0$ by (\ref{EQ2621}) we have  $u\in L_0 \cup L_1 \cup A_{<i_0} $ and,
by (\ref{EQ2631}),  we have $f(u)\not\in f[K]$. By (\ref{EQ2621}) and (\ref{EQ2630}), $u\in \bigcup _{k\leq n_0}L_k$ so
$f(u)\in f[\bigcup _{k\leq n_0}L_k]\setminus f[K]$ and, by (\ref{EQ2629}), $a_n \not\sim f(u)$, which, by (\ref{EQ2615}),
gives $f(a_n) \not\sim u$. Thus $f(a_n) \in R^{H_0}_{K_0}$.
\kdok
\noindent
Now, since $n\in S_{i_0,1}$ we have $f(a_n)\in f[A_1]$ and, by  (\ref{EQ2629}) and Subclaims 2.1 and 2.2,
$f(a_n)\in f[A_1] \cap R^{\{ p \} \cup H_0 \cup H_1}_{\;\;\emptyset\;\, \cup K_0 \cup K_1}$. Claim 2 is proved.
\kdok
\noindent
By Claims 1 and 2 (\ref{EQ2633}) is true and $A_0$ is a $(p,F,G)$-extendible copy below $A$.
\hfill $\Box$
\begin{lem}   \label{T2623}
For each maximal antichain $\A$ in $\P (R)$ and each finite set $F_0\subset R$
there is $S\in \P (R)$ containing $F_0$ and compatible with $\leq 2^{|F_0|}$ elements of ${\mathcal A}$.
\end{lem}

%==================================================================================================

\dok
We prove the lemma by induction on $|F_0|=k$. For $k=0$ this is trivial: take $S\in \A$.
Suppose that the statement is true for $k$. Let $|F_0|=k+1$, $p\in F_0$ and let
\begin{equation}\label{EQ2678}
F=F_0 \cap R^{ \{ p \} }_{\{ p \}} \;\;\mbox{ and }\;\; G= F_0 \cap  R^{ \{ p \} }_{\;\;\emptyset} .
\end{equation}
By Lemmas \ref{T2621} and  \ref{T2619} there is an isomorphism
$f:R ^{\{ p \}}_{\{ p \}}\rightarrow R ^{\{ p \}}_{\;\;\emptyset }$  satisfying (\ref{EQ2615}) and (\ref{EQ2666})
and there is a copy $B\in \P (R)$ satisfying $F\cup f^{-1}[G]\subset B\subset R^{ \{ p \} }_{\{ p \}}$ and such
that for each copy $A\in \P (R)$ satisfying
$F\cup f^{-1}[G] \subset A\subset B$ there are copies $A_0$ and $A_1$ satisfying (\ref{EQ2670}) and (\ref{EQ2633}).

\vspace{2mm}

\noindent
{\bf  Claim 1.}
$\D =\{ C\in \P (B) : \exists A' , A'' \in \A \;\; C\subset A' \cap f^{-1 }[A'' \cap R^{ \{ p \} }_{\;\;\emptyset} ]\}$ is a dense set
in the poset $\la \P (B ), \subset \ra$ (for each $E\in \P (B)$ there is $C\in \D$ such that $C\subset E$).

\vspace{2mm}

\noindent
{\it  Proof of Claim 1.}
Let $E\in \P (B)$. Since $E\in \P (R)$, by the maximality of $\A$ there are $A'\in \A$ and $C_1\in \P (R)$
such that $C_1\subset E\cap A'$. Since $f$ is an isomorphism we have $f[C_1]\in \P (R)$ and, again, there are $A'' \in \A$ and $C_2\in \P (R)$
such that $C_2\subset f[C_1]\cap A''$, which implies that for $C=f^{-1 }[C_2]$ we have $C\subset C_1 \subset E \subset B$ and, thus $C\in \P (B)$,
and $ C\subset f^{-1 }[A'' \cap R^{ \{ p \} }_{\;\;\emptyset} ]$. Since
$C\subset C_1 \subset A'$ we have $C\in \D$ and $C\subset C_1\subset E$.
\kdok
Let $\A ^*$ be a maximal antichain in the poset $\la \D , \subset \ra$.

\vspace{2mm}
\noindent
{\bf  Claim 2.}
$\A ^*$ is a maximal antichain in the poset $\la \P (B ), \subset \ra$.

\vspace{2mm}
\noindent
{\it  Proof of Claim 2.}
By the density of $\D$, $\A ^*$ is an antichain in $\la \P (B ), \subset \ra$. If $E\in \P (B)$, by Claim 1 there is $C\in \D$
such that $C\subset E$ and, by the maximality of $\A ^*$ in $\D$, there are $A\in \A ^*$ and $C_1\in \D$ satisfying $C_1 \subset C\cap A \subset E\cap A$.
Thus each $E\in \P (B)$ is compatible with some element of $\A^*$.
\kdok
Since $B\cong R $ (which implies $\P (B ) \cong  \P (R )$) and since $F\cup f^{-1}[G]\in [B]^k$ and $\A ^*$ is a maximal antichain in
$\P (B )$, by the induction hypothesis applied to $B$ there is  a set $A$ satisfying
\begin{equation}\label{EQ2683}
F\cup f^{-1}[G]\subset A \in \P (B)
\end{equation}
and compatible with $m\leq 2^k$ elements of $\A ^*$, say $C_1, \dots ,C_m$. Thus
\begin{equation}\label{EQ2684}
\forall C\in \A ^* \setminus \{ C_1 , \dots ,C_m \}\;\; A\perp C .
\end{equation}
Since $\A ^* \subset \D$, there are sets $A_1', A_1'' , \dots , A_m ',A_m ''\in \A $ such that
\begin{equation}\label{EQ2685}
\forall i\leq m \;\; C_i \subset A_i ' \cap f^{-1}[A_i '' \cap R^{ \{ p \} }_{\;\;\emptyset}] .
\end{equation}
By (\ref{EQ2683}) and Lemma \ref{T2619} there are sets $A_0$ and $A_1$ satisfying
\begin{equation}\label{EQ2681}
A_0\cup A_1\subset A \;\; \land  \;\; A_0 \cap A_1 =\emptyset  \;\; \land \;\;  F\subset A_0  \;\; \land  \;\; G\subset f[A_1] ,
\end{equation}
\begin{equation}\label{EQ2682}
S:= A_0 \cup \{ p \}\cup f[A_1]\in \P (R).
\end{equation}
By (\ref{EQ2681}) and (\ref{EQ2682}) we have $F_0 =F\cup \{ p \} \cup G  \subset S$ and it remains to be proved that
$S$ is compatible with $\leq 2m (\leq 2^{k+1})$-many elements of $\A$. Since $A_0 , A_1 \subset A$, by (\ref{EQ2684}) we have
\begin{equation}\label{EQ2686}
\forall C\in \A ^* \setminus \{ C_1 , \dots ,C_m \}\;\; (A_0 \perp C \land A_1 \perp C)
\end{equation}
and the proof will be finished when we show that
\begin{equation}\label{EQ2687}
\forall D\in \A \setminus \{ A_1', A_1'' , \dots , A_m ',A_m '' \}\;\; S \perp D.
\end{equation}
On the contrary, suppose that there are $D\in \A \setminus \{ A_1', A_1'' , \dots , A_m ',A_m '' \}$ and $C\in \P (R)$ such that
\begin{equation}\label{EQ2688}
C\subset S\cap D.
\end{equation}
By (\ref{EQ2682}), (\ref{EQ2688}) and since $R$ is strongly indivisible (see Fact \ref{T2600}(b)), at least one of the sets $C^0 =C\cap A_0$ and $C^1 =C\cap f[A_1]$ is a copy of $R$.

If $C^0 \in \P (R)$, then, since $C^0 \subset A_0 \subset A \subset B $ and $\A ^*$ is a maximal antichain in $\P (B)$,
there is $C^* \in \A ^*$ such that $C^* \not\perp C^0$, which implies $C^* \not\perp A_0$ thus, by  (\ref{EQ2686}),
$C^* =C_i$, for some $i\leq m$. By (\ref{EQ2685}) we have $C^* \subset A_i'$ and, since $C^0\subset C\subset D$ and $C^* \not\perp C^0$,
we have $A_i' \not\perp D$, which implies $D=A_i'$. But this contradicts our assumption concerning $D$.

If $C^1 \in \P (R)$, then, since $C^1 \subset f[ A_1]$, we have $f^{-1}[C^1] \subset A_1 \subset B $
and, since $f$ is an isomorphism, $f^{-1}[C^1] \in \P (B)$.
Since $\A ^*$ is a maximal antichain in $\P (B)$
there is $C^* \in \A ^*$ such that $C^* \not\perp f^{-1}[C^1]$ and, since  $f^{-1}[C^1] \subset A_1$, we have
$C^* \not\perp A_1$. Thus, by  (\ref{EQ2686}), $C^* =C_i$, for some $i\leq m$.
By (\ref{EQ2688}) we have $C^1\subset D$ and, since $C^1\subset R^{ \{ p \} }_{\;\;\emptyset}$, we have
$f^{-1}[C^1]\subset f^{-1}[D \cap R^{ \{ p \} }_{\;\;\emptyset}]$ which implies
\begin{equation}\label{EQ2689}
C_i =C^* \not\perp f^{-1}[D \cap R^{ \{ p \} }_{\;\;\emptyset}].
\end{equation}
By (\ref{EQ2685}), $C_i \subset f^{-1}[A_i'' \cap R^{ \{ p \} }_{\;\;\emptyset}]$ so, by (\ref{EQ2689}),
$f^{-1}[D \cap R^{ \{ p \} }_{\;\;\emptyset}] \not\perp f^{-1}[A_i'' \cap R^{ \{ p \} }_{\;\;\emptyset}]$
and, hence, there is $E\in \P (R)$ such that
$E \subset f^{-1}[D \cap R^{ \{ p \} }_{\;\;\emptyset}] \cap f^{-1}[A_i'' \cap R^{ \{ p \} }_{\;\;\emptyset}]$.
But then $\P (R)\ni f[E] \subset D \cap A_i''$
and we have $D \not\perp A_i''$, which implies $D=A_i''$. This is a contradiction. Thus (\ref{EQ2687}) is true and the proof is finished.
\kdok
\noindent
{\bf Proof of Theorem \ref{T2635}.}
Let $F_0\in [R]^{<\o}$ and let $\A$ be a maximal antichain in $\P (R)$.
First we prove that
$$\D =\{ C\in \P (R): \exists A\in \A \;\; \exists H\subset F_0 \;\; C\subset A \cap R^{F_0}_H \}$$
is a dense set in the poset $\P (R)$. If $B\in \P (R)$, then, by Fact \ref{T2600}(a) and (c), $B\setminus F_0 \in \P (R)$ and
$B\setminus F_0 = \bigcup _{H\subset F_0}B\cap R^{F_0}_H $
so, since
$B\setminus F_0$ is strongly indivisible (Fact \ref{T2600}(b)), there is $H_0 \subset F_0$ such that $B\cap R^{F_0}_{H_0} \in \P (R)$. By the maximality of
$A$ there are $A_0 \in \A$ and $C\in \P (R)$ such that $C\subset B\cap R^{F_0}_{H_0}\cap A_0$. Thus $C\in \D$ and $C\subset B$.

Let $\A ^*$ be a maximal antichain in the poset $\la \D ,\subset \ra$. Clearly $\A ^*$ is a maximal antichain in the poset $\la \P (R) ,\subset \ra$
and, by Lemma \ref{T2623}, there is $S\in \P (R)$ containing $F_0$ and compatible with $m\leq 2^{|F_0|}$ elements of ${\mathcal A}^*$,
say $C_1 , \dots , C_m$. Next we prove that
\begin{equation}\label{EQ2693}
\forall H \subset F_0 \;\; \exists _1 i \leq m \;\; (C_i \subset R^{F_0}_H \land C_i \not\perp S^{F_0}_H ).
\end{equation}
Let $H \subset F_0$. Since $S^{F_0}_H\in \P (R)$ there is $C\in \A ^*$ such that $C\not\perp S^{F_0}_H$, which implies $C\not\perp S$ and, hence,
$C=C_{i_H}$, for some $i_H\leq m$. Since $C_{i_H} \in \D$ there is $H'\subset F_0$ such that $C_{i_H}\subset R^{F_0}_{H'} $ and, since $C_{i_H}$ is compatible with
$S^{F_0}_H \subset R^{F_0}_H$, we have $R^{F_0}_H \not\perp R^{F_0}_{H'}$, which implies $H'=H$ and $C_{i_H}\subset R^{F_0}_H$. Thus, since $m\leq 2^{|F_0|}$,
$H \mapsto i_H$ is an bijection from $P(F_0)$ to $\{ 1,2, \dots , m \}$ and (\ref{EQ2693}) is true.

Now we prove
\begin{equation}\label{EQ2694}\textstyle
S_1 =F_0 \cup \bigcup _{H\subset F_0}(C_{i_H}\cap S^{F_0}_H) \in \P (R).
\end{equation}
Suppose that $S_1\not\in \P (R)$. Then
$S\setminus S_1 =\bigcup _{H\subset F_0}(S^{F_0}_H\setminus C_{i_H}) \in \P (R)$ ($S$ is strongly indivisible) and, by the maximality of $\A ^*$,
there are $C^*\in \A ^*$ and $C\in \P (R)$ such that $C\subset C^* \cap (S\setminus S_1)=\bigcup _{H\subset F_0}C^* \cap (S^{F_0}_H\setminus C_{i_H})$.
Thus, since $C$ is  strongly indivisible, there is $H_0\subset F_0$ such that $C_1 =C\cap C^* \cap (S^{F_0}_{H_0}\setminus C_{i_{H_0}})\in \P (R)$,
which implies that $S$ is compatible with $C^*\in \A ^* \setminus \{ C_{i_H} : H\subset F_0 \}$. But, by (\ref{EQ2693}),
$\{ C_{i_H} : H\subset F_0 \}= \{ C_i :i\leq m \}$, a contradiction. Thus (\ref{EQ2694}) is true.

Finally we prove (\ref{EQ2692}). For $H\subset F_0$ we have $C_{i_H} \in \A ^* \subset \D$ and, hence, there is $A\in \A$ such that $C_{i_H}\subset A$. Thus
by (\ref{EQ2694}),
$(S_1)^{F_0}_H =C_{i_H}\cap S^{F_0}_H\subset A$.
\hfill $\Box$
\section{Fusion for $\P (R)$}\label{S5}
If $\la R, \sim \ra$ is the Rado graph and $\D= \la \D _n :n\in \o \ra$ a sequence of subsets of $\P (R)$, then  a copy $L\in \P (R)$ will be called
{\it a fusion of $\D$} if and only if there exists a labeling $\la \{ L_n : n\in \o \}, \{ q ^{\bigcup _{i<n}L_i}_K : n\in \o \land K\subset \bigcup _{i<n}L_i\}\ra$
of $L$ such that
\begin{equation}\label{EQ2602}\textstyle
\forall n\in \o \;\; \forall  K \subset \bigcup _{i<n}L_i\;\;  \exists D \in \D _n \;\;L^{\bigcup _{i<n}L_i}_K \subset D .
\end{equation}
\begin{te}   \label{T2637}
If  $\D= \la \D _n \!:n\in \o \ra$ is a sequence subsets of $\P (R)$ which are dense below $A\in \P (R)$,
then the set ${\mathcal F}=\{L : L \mbox{ is a fusion of }\D \}$ is dense below $A$.
\end{te}
\dok
Let $B\in \P (R)$ and $B\subset A$. In order to construct an $L\in {\mathcal F}\cap \P (B)$ by
recursion we define a sequence $\la \la S_n ,L_n \ra : n\in \o \ra$
such that for each $n \in \o$

(i) $S_n \in \P (B)$,

(ii) $S_{n+1}\subset S_n$,

(iii) $\bigcup _{i\leq n}L_i\subset S_n$,

(iv) $L_n =\{ q ^{\bigcup _{i<n}L_i}_K : K\subset \bigcup _{i<n}L_i \}$, where $q ^{\bigcup _{i<n}L_i}_K \in  (S_n) ^{\bigcup _{i<n}L_i}_K$, and

(v) $\forall K\subset \bigcup _{i<n}L_i\;\; \exists D \in \D _n \;\;(S_n) ^{\bigcup _{i<n}L_i}_K\subset D$.

\noindent
Since $B\subset A$ the set $\D _0$ is dense below $B$. We choose $S_0\in \D _0 $ such that $S_0\subset B$, take $q^\emptyset _\emptyset \in S_0$,
define $L_0= \{ q^\emptyset _\emptyset\}$ and conditions (i) - (v) are satisfied.

Suppose that a sequence $\la\la S_i, L_i \ra :i<n\ra $ satisfies conditions (i) - (v). Then $S_{n-1}\in \P (B)$ and, hence,
the set $\D _n '=\{ D\in \D _n : D\subset S_{n-1}) \}$
is dense below $S_{n-1}$. Let $\A _n$ be a maximal antichain in $\la \D _n ', \subset\ra$. Clearly $\A _n$ is a maximal antichain in
the poset $\la \P (S _{n-1}) , \subset\ra$ and, by (iii), $\bigcup _{i< n}L_i\subset S_{n-1}$ so, by Theorem \ref{T2635} applied to $S_{n-1}$,
there is a set $S_n$ satisfying
\begin{equation}\label{EQ2639}\textstyle
\bigcup _{i< n}L_i\subset S_n \in \P (S_{n-1}), \mbox{ and}
\end{equation}
\begin{equation}\label{EQ2640}\textstyle
\forall K\subset \bigcup _{i< n}L_i \;\; \exists D\in \A _n \;\;(S_n) ^{\bigcup _{i<n}L_i}_K \subset D .
\end{equation}
By (\ref{EQ2639}) conditions (i) and (ii) are satisfied and, since $S_n \cong R$, for $K\subset \bigcup _{i<n}L_i$ we choose
$q ^{\bigcup _{i<n}L_i}_K \in  (S_n) ^{\bigcup _{i<n}L_i}_K$ and define
$L_n =\{ q ^{\bigcup _{i<n}L_i}_K : K\subset \bigcup _{i<n}L_i \}$, so (iii) and (iv) are satisfied too. By (\ref{EQ2640}),
for $K\subset \bigcup _{i< n}L_i$ there is $D\in \A _n \subset \D _n $ such that $(S_n) ^{\bigcup _{i<n}L_i}_K \subset D$.
Thus (v) is true and the recursion works.

We show that $L:= \bigcup _{n\in \o}L_n\in {\mathcal F}$.
By (i) we have $S_n \subset A$ and, by (iv), for $K\subset \bigcup _{i<n}L_i$ we have
$q ^{\bigcup _{i<n}L_i}_K \in  (S_n) ^{\bigcup _{i<n}L_i}_K \subset A ^{\bigcup _{i<n}L_i}_K$ and, by Lemma \ref{T2633}(b),
$\la \{ L_n : n\in \o \}, \{ q ^{\bigcup _{i<n}L_i}_K : n\in \o \land K\subset \bigcup _{i<n}L_i\}\ra$ is a labeling of $L$.
By (ii) and (iii), for $n\in \o$ and $K\subset \bigcup _{i<n}L_i$ we have $L\subset S_n$, which together with (v) implies
that there is $D\in \D _n$ such that $L ^{\bigcup _{i<n}L_i}_K \subset  (S_n) ^{\bigcup _{i<n}L_i}_K \subset D$ and (\ref{EQ2602}) is true as well.
Thus $L\in {\mathcal F}$ and, by (i) and (ii), $L\subset B$;
so, ${\mathcal F}$ is dense below $A$.
\kdok
The following statement is an improvement of Theorem \ref{T2635}.
\begin{cor}\label{T2638}
For each sequence $\A=\{ \A _n : n\in \o \}$ of maximal antichains in the poset $\P (R)$ there is a maximal antichain $\A '$ in $\P (R)$
consisting of fusions of $\A$.
\end{cor}
\dok
By the assumption, the sets $\D _n =\{ D\subset R : \exists A\in \A _n \; D\subset A\}$, $n\in \o$, are dense in $\P (R)$ and, by Theorem \ref{T2637},
the corresponding set of fusions ${\mathcal F}$ is a dense set as well. If $\A '$ is a maximal chain in ${\mathcal F}$ it is a maximal chain in $\P (R)$.
\hfill $\Box$
\section{New reals and a factorization}\label{S6}
Clearly, the $\P (R)$-name $\r =\{ \la \check{p}, R^{\{ p \}}_{\{ p \}}\ra : p\in R \}$ is a name for a subset of $R$ and, since $|R|=\o$, $\r$
can be regarded as a name for a real.
\begin{te}   \label{T2631}
The name $\r$ codes a new real, that is, $R\Vdash \r \in P(\check{R}) \setminus V$.
\end{te}
\dok
Let $G$ be a $\P (R)$-generic filter over $V$.
Suppose that $\r _G =S$ for some $S\in P(R)\cap V$. Then $A\Vdash \r =\check{S}$, for some $A\in G$,
which implies that $A\Vdash \check{p}\in \r$, for all $p\in S$, and  $A\Vdash \check{p}\not\in \r$, for all $p\in R\setminus S$.
Since $A\Vdash \check{p}\in \r$ iff $A\cap R^{\{ p\}}_{\;\;\emptyset }\in \I _R$
and $A\Vdash \check{p}\not\in \r$ iff $A\cap R^{\{ p\}}_{\{ p \} }\in \I _R$ we have
\begin{equation}\label{EQ2661}\textstyle
\forall p\in S \;\; (A\cap R^{\{ p\}}_{\;\;\emptyset }\in \I _R ) \;\;\land\;\;
\forall p\in R\setminus S \;\; (A\cap R^{\{ p\}}_{\{ p \} }\in \I _R)
\end{equation}
Let $p\in A$. If $p\in S$, then, since $p\in A\in \P (R)$, we have $\P (R)\ni A^{\{ p\}}_{\;\;\emptyset }=A\cap R^{\{ p\}}_{\;\;\emptyset }$,
which is impossible by (\ref{EQ2661}). If $p\in R\setminus S$, then, $\P (R)\ni A^{\{ p\}}_{\{ p \} }=A\cap R^{\{ p\}}_{\{ p\} }$,
which is impossible by (\ref{EQ2661}). A contradiction.
\kdok
\noindent
If $G$ is a $\P (R)$-generic filter over the ground model $V$ (of ZFC),
then, by Theorem \ref{T2631}, $\r _G \not\in V$ and (see \cite{Jech}, p.\ 265)
there is a forcing $\P$ and a $\P $-generic filter over $V$, $H$, such that $V[\r _G]=V_\P [H]$. Thus (see \cite{Jech1}, p.\ 48) there is a $\P$-name for a poset
$\pi$ such that the generic extension $V_{\P (R)}[G]$ is equal to the two-step extension
$(V_\P [H])_{\pi _H}[H_1]=(V[\r _G])_{\pi _H}[H_1]$, where $H_1$ is a $\pi _H$-generic filter over $V_\P [H]$. In the sequel we show that $\pi$ is a name
for an $\o$-distributive forcing.
\begin{te}   \label{T2634}
Let $\k $ be an infinite cardinal and $G$ a $\P (R)$-generic filter over the ground model $V$.
If $x\in V_{\P (R)}[G]$, where $x:\o \rightarrow \k$,
then $x\in V[\r _G]$.
\end{te}
\dok
Let $\t$ be a $\P (R)$-name such that $x=\t _G$. Then there is $A\in G$ such that $A\Vdash \t : \check{\o }\rightarrow \check{\k }$ and
first we prove that
\begin{equation}\label{EQ2690}\textstyle
\forall B \in \P (A)\;\; \exists L \in \P (B) \;\; L\Vdash \t \in V[\r ].
\end{equation}
Let $B \in \P (A)$. Since $A\Vdash  \forall n\in \check{\o } \; \exists \a \in \check{\k }\; \t (\check{n})=\check{\a }$, for each $n\in \o$ we have:
for each $C\in \P (A)$ there are $D\in \P (C)$ and $\a \in \k$ such that $D\Vdash \t (\check{n})=\check{\a }$. This means that the sets
$\D _n :=\{ D\in \P (A): \exists \a \in \k \;\; D\Vdash \t (\check{n})=\check{\a }\}$, $n\in \o$, are dense below $A$. By Theorem
\ref{T2637}, the set ${\mathcal F}$ of fusions is dense below $A$ and, hence, there is $L=\bigcup _{n\in \o }L_n \in {\mathcal F}$ such that $L\in \P (B)$.
By (\ref{EQ2602}), for
$n\in \o$ and $K \subset \bigcup _{i<n}L_i$ there is $D \in \D _n $ such that $L^{\bigcup _{i<n}L_i}_K \subset D $ and, hence, there is (clearly unique)
$\a \in \k$ such that $L^{\bigcup _{i<n}L_i}_K \Vdash \t (\check{n})=\check{\a }$.
Thus we obtain a family of ordinals $\{ \a^{\bigcup _{i<n}L_i}_K :n\in \o \land K\subset \bigcup _{i<n}L_i\}$ indexed by elements $q^{\bigcup _{i<n}L_i}_K$
of $L$ such that
\begin{equation}\label{EQ2662}\textstyle
L^{\bigcup _{i<n}L_i}_K\Vdash \t (\check{n})=\check{\a^{\bigcup _{i<n}L_i}_K}.
\end{equation}
In order to prove that $L \Vdash \t \in V[\r ]$ we assume that $H$ is a $\P (R)$-generic filter over $V$ containing $L$ and
we reconstruct $\t _H$ inside $V[\r _H]$ showing that for each $n\in \o$
$$
\t _H (n)= \a ^{\bigcup _{i<n}L_i} _{  (\bigcup _{i<n}L_i ) \cap \r _H},
$$
which will, by (\ref{EQ2662}), follow from  $L^{\bigcup _{i<n}L_i}_{(\bigcup _{i<n}L_i ) \cap \r _H}\in H$.
Clearly $R^{\{ p \}}_{\{ p \}}\in H $, for each $p\in  (\bigcup _{i<n}L_i) \cap \r _H $, and
$R^{\{ p \}}_{\;\;\emptyset }\in H $, for each $p\in (\bigcup _{i<n}L_i )\setminus \r _H  $.
Thus
$$\textstyle
R^{\bigcup _{i<n}L_i}_{(\bigcup _{i<n}L_i ) \cap \r _H}
=\bigcap _{p\in  (\bigcup _{i<n}L_i) \cap \r _H}R^{\{ p \}}_{\{ p \}}
\cap
\bigcap _{p\in (\bigcup _{i<n}L_i )\setminus \r _H }R^{\{ p \}}_{\;\;\emptyset }\in H$$
and, since $L \in H$, we have
$L^{\bigcup _{i<n}L_i}_{(\bigcup _{i<n}L_i ) \cap \r _H}  =L \cap  R^{\bigcup _{i<n}L_i}_{(\bigcup _{i<n}L_i ) \cap \r _H} \in H$.
So $\t _H \in V[\r _H]$ and we proved that $L \Vdash \t \in V[\r ]$, which completes the proof of (\ref{EQ2690}).

Now, since $A\in G$, by (\ref{EQ2690}) there is $L\in G$ satisfying $L\Vdash \t \in V[\r ]$ and, hence,
$x= \t _G \in V[\r _G]$.
\hfill $\Box$
\section{The $\aleph _0$-covering and the Sacks property}\label{S7}
\noindent
For a cardinal $\k \geq \o$ and a sequence of positive integers $\la k_n : n\in \o \ra$
a mapping $s:\o \rightarrow [\k ]^{<\o }$ will be called an {\it $\la k_n \ra$-slalom in $\k$} iff $|s(n)|\leq k_n$, for each $n\in \o$.
$\Sl _{\la k_n \ra}(\k )$ will denote the set of all such mappings.

A pre-order $\P$ has the {\it Sacks property} iff
there is a sequence $\la k_n \ra \in \N ^\o$ such that in each generic extension $V_\P [G]$ for each $x : \o \rightarrow \o$ there is
$s\in V \cap \Sl _{\la k_n \ra}(\o )$ (or, equivalently, $s\in V \cap \Sl _{\la 2^n \ra}(\o )$) such that $x(n)\in s(n)$, for each $n\in \o$.

A pre-order $\P$ has the {\it $\aleph _0$-covering property} iff in each generic extension $V_\P [G]$ each countable set $X$ of ordinals is contained in
a countable set $A\in V$.

We note that the Cohen forcing has the $\aleph _0$-covering property (it is a ccc poset) but does not have the Sacks property, while, under CH, the Namba forcing
has the Sacks property (since it does not produce new reals) but does not have the $\aleph _0$-covering (since it adds a cofinal mapping from $\o$ to $\o _2$, see \cite{Jech});
the Sacks forcing has both of these properties and we show that the same holds
for the forcing $\P (R)$.

We recall that a complete Boolean algebra $\B$ is weakly distributive
 iff for each cardinal $\k$ in each generic extension $V_{{\mathbb B}}[G]$ for each function $x:\o \rightarrow \k$ there is a slalom
$s: \o \rightarrow [\k ]^{<\o }$ belonging to $V$ and such that $x(n)\in s(n)$, for all $n\in \o$.
\begin{te}   \label{T2625}
(a) If $\k $ is an infinite cardinal and $G$ a $\P (R)$-generic filter over the ground model $V$,
then for  each function
$x : \o \rightarrow \k$ belonging to $V_{\P (R)}[G]$ there exists a slalom
$s\in V \cap \Sl _{\la m_n \ra}(\k )$ such that $x(n)\in s(n)$, for each $n\in \o$.

(b) The forcing $\P (R)$ has the $\aleph _0$-covering property and, hence,  preserves $\o _1$.

(c) The forcing $\P (R)$ has the Sacks property.

(d) The algebra $\ro \sq \P (R)$ is a weakly distributive complete Boolean algebra.
\end{te}
\dok
(a) We have to prove that for each $\P (R)$-name $\t$
\begin{equation}\label{EQ2695}\textstyle
R \Vdash \t : \check{\o }\rightarrow \check{\k } \Rightarrow
\exists s\in (( \Sl _{\la m_n \ra}(\check{\k }))^V)\check{\;}\;\; \forall n\in \check{\o } \;\;\t (n)\in s(n).
\end{equation}
Thus, working in $V$ we show that for each $A \in \P (R)$ satisfying $A\Vdash \t : \check{\o }\rightarrow \check{\k }$ there are $L\in \P (A)$ and
$s\in \Sl _{\la m_n \ra}(\k )$ such that
\begin{equation}\label{EQ2696}\textstyle
\forall n\in \o \;\; L \Vdash   \t (\check{n})\in \check{s}(\check{n}).
\end{equation}
First, exactly as in the proof of Theorem \ref{T2634} we find $L=\bigcup _{n\in \o }L_n \in \P (A)$ such that
\begin{equation}\label{EQ2699}\textstyle
\forall n\in \o \;\; \forall K\subset \bigcup _{i<n}L_i\;\; \exists _1 \a^n_K\in \k \;\; L^{\bigcup _{i<n}L_i}_K\Vdash \t (\check{n})=\check{\a^n_K}.
\end{equation}
Let the mapping $s:\o \rightarrow [\k ]^{<\o }$ be defined by $s(n)=\{ \a ^n_K :K\subset \bigcup _{i<n}L_i \}$. Since
$\{ L^{\bigcup _{i<n}L_i}_K : K\subset \bigcup _{i<n}L_i\} $ is a maximal antichain below $L$ in $\P (R)$, by (\ref{EQ2699})
we have $L\Vdash \t (\check{n})\in \check{s}(\check{n})$ and (\ref{EQ2696}) is true. Finally,
$s\in \Sl _{\la m_n \ra}(\k )$ because $|L_i|=m_i$ and $|s(n)|\leq |P(\bigcup _{i<n}L_i)| = 2^{\sum _{i<n}|L_i|}=2^{\sum _{i<n}m_i}=m_n$.

(b) If $X\in V_{\P (R)}[G]\cap [\kappa ]^\o$ and $x:\o \rightarrow X$ is a bijection, then by (a) we have $X=x[\o ]\subset \bigcup _{n\in \o}s(n)\in V$,
because $s\in V$.

(c) By (a) each function $x : \o \rightarrow \o$ is contained in an $s\in V \cap \Sl _{\la m_n \ra}(\o )$.
\hfill $\Box$
\section{Tree-ordered copies of the Rado graph}\label{S8}
Here we show that each labeling of a copy $L$ of the Rado graph $\la R, \sim \ra$ induces a reversed tree order on $L$ in a natural way.
This construction will be used in the next section.
So, let  $\CL=\la \{ L_n :n\in \o \}, q \ra$ be a labeling of $L$, that is
\begin{equation}\label{EQ2659}\textstyle
L_n =\{ q ^{\bigcup _{i<n}L_i}_K : K\subset \bigcup _{i<n}L_i \}, \mbox{ where } q ^{\bigcup _{i<n}L_i}_K \in  L^{\bigcup _{i<n}L_i}_K .
\end{equation}
Using the labeling $\CL$ we define the binary relation $\leq _{L,\CL}  $ (we will write shortly $\leq _L$) on $L$ by:
\begin{equation}\label{EQ2622}\textstyle
q ^{\bigcup _{i<m}L_i}_{K '} \leq _L  q ^{\bigcup _{i<n}L_i}_{K''} \;\;\Leftrightarrow \;\;L^{\bigcup _{i<m}L_i}_{K '} \subset  L ^{\bigcup _{i<n}L_i}_{K''}.
\end{equation}
Since $L^\emptyset _\emptyset =L$ we have $q\leq _L q^\emptyset _\emptyset$, for all $q\in L$ and, clearly,
the relation $\leq _L$ is reflexive, transitive and, by Lemma \ref{T2620}(c) antisymmetric. Thus, $\la L, \leq _L \ra$ is a partial order
with the largest element $q^\emptyset _\emptyset$.
Lemma \ref{T2620}(d) applied to $L$ gives $L^{\bigcup _{i<m}L_i}_{K '} \subset  L ^{\bigcup _{i<n}L_i}_{K''}$ if and only if
$\bigcup _{i<m}L_i \supset \bigcup _{i<n}L_i$, $K'\supset K''$ and $K' \cap \bigcup _{i<n}L_i =K''$ if and only if $m\geq n $ and $K' \cap \bigcup _{i<n}L_i =K''$.
Thus, by (\ref{EQ2622}) we have
\begin{equation}\label{EQ2605}\textstyle
q ^{\bigcup _{i<m}L_i}_{K '} \leq _L  q ^{\bigcup _{i<n}L_i}_{K''} \;\;\Leftrightarrow \;\; m\geq n \;\;\land \;\; K' \cap \bigcup _{i<n}L_i =K'' .
\end{equation}
In order to state the
following theorem we introduce a convenient notation. For $p,q\in L$ let $p\prec _L q$ denote that $p$ is an immediate predecessor of $q$ in $\la L , \leq _L  \ra$ and let
$$
\Ip _{\la L , \leq _L  \ra}(q)=\{ p\in L: p\prec _L q \}.
$$
For $q ^{\bigcup _{i<n}L_i}_K\in L_n$ and $K_1\subset L_n$ let
$q ^{\bigcup _{i<n}L_i}_K \Rsh ^{L_n}_{K_1}$ denote the element $q ^{\bigcup _{i< n+1}L_i}_{K\cup K_1}$ of $L_{n+1}$.
For simplicity, the intervals $(p,q]_{\la L, \leq _L  \ra}$ will be denoted by $(p,q]$.
\begin{te}   \label{T2627}
For each $n\in\o$ and $K \subset \bigcup _{i<n}L_i$ in the poset $\la L, \leq _L  \ra$ we have:

(a) $(q ^{\bigcup _{i<n}L_i}_{K },q^\emptyset _\emptyset ]= \{ q ^{\bigcup _{i<m}L_i}_{K  \cap \bigcup _{i<m}L_i}: m<n \} $;

(b) $\la L, \leq _L  \ra$ is a reversed tree with the top $q^\emptyset _\emptyset $ and the set $L_n$ is its $n$-th level;

(c) $(-\infty , q ^{\bigcup _{i<n}L_i}_K ]= L^{\bigcup _{i<n}L_i}_K $;

(d) $\Ip _{\la L , \leq _L  \ra}(q ^{\bigcup _{i<n}L_i}_K) =\{ q ^{\bigcup _{i<n}L_i}_K \Rsh ^{L_n}_{K_1} : K_1 \subset L_n \}$;

(e) $\la L, \leq _L  \ra$ is a finitely branching reversed tree without minimal nodes. In fact each element of $L_n$ has
$2^{|L_n|}=2^{m_n}$ immediate predecessors.
\end{te}
\dok
(a)
If $q ^{\bigcup _{i<n}L_i}_{K} < _L  q ^{\bigcup _{i<m }L_i}_{K'}  $, then $n\geq m$ and $K'= K \cap \bigcup _{i<m}L_i $.
Since $q ^{\bigcup _{i<n}L_i}_{K} \neq _L  q ^{\bigcup _{i<m }L_i}_{K'} $, we have  $m<n$. Thus ``$\subset$" is true and
``$\supset$" is obvious.

(b) If $m_1 < m_2 <n$, then
$q ^{\bigcup _{i<m_2}L_i}_{K  \cap \bigcup _{i<m_2}L_i} <_L q ^{\bigcup _{i<m_1}L_i}_{K  \cap \bigcup _{i<m_1}L_i}$ and by (a) the interval
$(q ^{\bigcup _{i<n}L_i}_{K },q^\emptyset _\emptyset ]$ is a chain of size $n$. Thus
$\la L, \leq _L  \ra$ is a reversed tree, $\height (q ^{\bigcup _{i<n}L_i}_{K})$ $=n$
and, hence, $L_n =\{ q ^{\bigcup _{i<n}L_i}_{K} : K\subset \bigcup _{i<n}L_i \}$ is the $n$-th level of $\la L, \leq _L  \ra$.

(c) If  $q ^{\bigcup _{i<m}L_i}_{K'}\leq _L q ^{\bigcup _{i<n}L_i}_K$, then,  by (\ref{EQ2659}) and (\ref{EQ2622}),
$q ^{\bigcup _{i<m}L_i}_{K'} \in L ^{\bigcup _{i<m}L_i}_{K'}\subset L ^{\bigcup _{i<n}L_i}_K$.
Conversely, if $q ^{\bigcup _{i<m}L_i}_{K'}\in L ^{\bigcup _{i<n}L_i}_K$, then $q ^{\bigcup _{i<m}L_i}_{K'}\in L_m \setminus \bigcup _{i<n}L_i$, which implies
$m\geq n$. Since $q ^{\bigcup _{i<m}L_i}_{K'}\in R^{\bigcup _{i<m}L_i}_{K'}\cap R ^{\bigcup _{i<n}L_i}_K $, by Lemma \ref{T2620} we have
$K' \cap \bigcup _{i<n}L_i = K \cap \bigcup _{i<m}L_i$ and, since $K\subset \bigcup _{i<n}L_i\subset \bigcup _{i<m}L_i$ we obtain $K' \cap \bigcup _{i<n}L_i = K$.
Thus,  by (\ref{EQ2605}), $q ^{\bigcup _{i<m}L_i}_{K'}\leq _L q ^{\bigcup _{i<n}L_i}_K$.

(d) If $q ^{\bigcup _{i<p}L_i}_F \prec _L   q ^{\bigcup _{i<n}L_i}_K$, then by (\ref{EQ2605}) $p\geq n$ and $F\cap\bigcup _{i<n}L_i =K$.
Since $p=n$ would imply the equality, we have $p\geq n+1$. For $p>n+1$ we would have
$    q ^{\bigcup _{i<p}L_i}_F
<_L  q ^{\bigcup _{i<n+1}L_i}_{F\cap \bigcup _{i<n+1}L_i}
<_L  q ^{\bigcup _{i<n}L_i}_{F\cap \bigcup _{i<n}L_i}
=    q ^{\bigcup _{i<n}L_i}_K$
which is not true. Thus $p=n+1$ and
$q ^{\bigcup _{i<p}L_i}_F= q ^{\bigcup _{i<n+1}L_i}_F=q ^{\bigcup _{i<n}L_i \cup L_n}_{K \cup (F\cap L_n)}=q ^{\bigcup _{i<n}L_i }_{K }\Rsh ^{L_n}_{F\cap L_n}$.

Conversely, for $K_1 \subset L_n$ we have $(K\cup K_1 )\cap \bigcup _{i< n}L_i=K$ and, by (\ref{EQ2605}),
$q ^{\bigcup _{i<n}L_i}_K \Rsh ^{L_n}_{K_1}= q ^{\bigcup _{i< n+1}L_i}_{K\cup K_1}\leq _L q ^{\bigcup _{i<n}L_i}_K$.
If $q ^{\bigcup _{i< n+1}L_i}_{K\cup K_1}\leq _L  q ^{\bigcup _{i<p}L_i}_F \leq _L   q ^{\bigcup _{i<n}L_i}_K$,
then $n+1 \geq p \geq n$.
So, if $p=n+1$, then $K\cup K_1 =(K\cup K_1 )\cap \bigcup _{i< n+1}L_i=F$, which implies
$ q ^{\bigcup _{i<p}L_i}_F= q ^{\bigcup _{i< n+1}L_i}_{K\cup K_1}$.
If $p=n$, then $K= F\cap \bigcup _{i<n}L_i =F\cap \bigcup _{i<p}L_i=F$,
thus $q ^{\bigcup _{i<p}L_i}_F =   q ^{\bigcup _{i<n}L_i}_K$. So
$(q ^{\bigcup _{i< n+1}L_i}_{K\cup K_1}, q ^{\bigcup _{i<n}L_i}_K)=\emptyset$, that is
$q ^{\bigcup _{i<n}L_i}_K \Rsh ^{L_n}_{K_1}\prec_L q ^{\bigcup _{i<n}L_i}_K$.
Clearly (e) follows from (d).
\kdok
Thus, by Lemma \ref{T2639} each copy of $R$ has infinitely many labelings and the corresponding induced reversed tree orderings.
By Theorem \ref{T2627}(c) and Corollary \ref{T2638}, we have
\begin{cor}\label{T2640}
For each sequence $\A=\{ \A _n : n\in \o \}$ of maximal antichains in $\P (R)$ there is a maximal antichain $\A '$ in $\P (R)$
such that each $L\in \A '$ has a labeling such that for each $n\in \o$ the set $\{ (-\infty , q]: q\in \Lev _n \la L,\leq _L\ra \}$ refines $\A _n\upharpoonright L$.
%the maximal antichain $\A _n \upharpoonright L$ in $\P (L)$.
\end{cor}
\section{Strong subtrees of the ordered Rado graph are large}\label{S9}
Let $T=\bigcup _{n\in \o} \Lev _n (T)$ be a tree of height $\o $. A subset $\CS $ of $T$ is called a {\it strong subtree}
of $T$ iff

(sst1) $\CS$ has the unique root,

(sst2) There is a set $S=\{ n_k : k\in \o \}\in [\o ]^\o$ such that $\emptyset \neq \Lev _k (\CS ) \subset \Lev _{n_k }(T)$, for each $k\in \o$,
($S$ is called the {\it level set of} $\CS$),

(sst3) If $s\in \Lev _k (\CS )$, then for each $T$-immediate successor $t$ of $s$ there is a unique $s_t \in \Lev _{k+1} (\CS )$ such that $s_t\geq t$.

We will use the following consequence of the Halpern-L\"{a}uchli Theorem (see \cite{Hal}, \cite{Tod}).
\begin{fac}[Halpern-L\"{a}uchli]  \label{T2628}
If $T$ is a countable finitely branching tree with one root and without maximal nodes,
then for each finite coloring of $T$ there is a monochromatic strong subtree of $T$.
\end{fac}
By (b) and (e) of Theorem \ref{T2627} and the obvious dual of Fact \ref{T2628} for reversed trees we have
\begin{te}   \label{T2636}
If $\la R, \sim \ra$ is the Rado graph, $L\subset R$ a copy of $R$, $\CL$ a labeling of $L$ and $\leq _L$ the corresponding order on $L$,
then for each finite coloring of the set  $L$ there is a monochromatic strong reversed subtree $\CS$ of the reversed tree $\la L,\leq _L  \ra$.
\end{te}
Now we show that each strong reversed subtree of $\la L,\leq _L  \ra$ contains a copy of $R$.
\begin{te}   \label{T2629}
Let $\la R, \sim \ra$ be the Rado graph, $L\in \P ( R)$ a copy of $R$ inside $R$,
$$\textstyle
\CL _L= \Big\la \{ L _n : n\in \o \} , \{ q ^{\bigcup _{i<n}L_i}_K : n\in \o \land K\subset \bigcup _{i<n}\L _i\} \Big\ra
$$
a labeling of $L$, $\leq _L$ the corresponding reversed tree order on $L$ and
$\CS$ a strong reversed subtree of $\la L,\leq _L  \ra$.
Then there is a copy $\L \in \P (R)$ satisfying $\L \subset \CS$ and there is a labeling of $\L$
$$\textstyle
\CL _\L= \Big\la \{ \L _k : k\in \o \} , \{ p ^{\bigcup _{j<k}\L _j}_H : k\in \o \land H\subset \bigcup _{j<k}\L _j\} \Big\ra
$$
such that the orders $\leq _\L$ and $\leq _L$ coincide on $\L$ and for $k\in \o$ and $H\subset \bigcup _{j<k}\L _j$,
\begin{equation}\label{EQ2613}\textstyle
\L ^{\bigcup _{j<k}\L _j}_H
=(-\infty , p ^{\bigcup _{j<k}\L _j}_H ]_{\la \L , \leq _\L  \ra}
=\L \cap (-\infty , p ^{\bigcup _{j<k}\L _j}_H ]_{\la L , \leq _L  \ra} .
\end{equation}
\end{te}
\dok
Let $S=\{ n_k:k\in \o \}$ be the level set of $\CS$, where $n_0<n_1<n_2\dots$.
For $q\in \Lev_k (\CS )=\CS \cap L_{n_k}$ and $p\in \Ip _{\la L , \leq _L  \ra}(q)$ let $\pi _\CS (p)$ denote the unique element of
$\Lev_{k+1} (\CS )=\CS \cap L_{n_{k+1}}$
satisfying $\pi _\CS (p) \leq _L p$ (such an element exists uniquely by the dual of (sst3) for reversed trees).

By recursion for $k\in \o$ we define $\L _k\subset R$ and $p^{\bigcup _{j<k}\L _j}_H \in R$, $H\subset \bigcup _{j<k}\L _j$, such that

($\L 1$) $\L _k \subset \CS \cap L_{n_k}$,

($\L 2$) $\L _k =\{ p^{\bigcup _{j<k}\L _j}_H : H\subset \bigcup _{j<k}\L _j\}$,

($\L 3$) $p^{\bigcup _{j<k}\L _j}_H =\pi _\CS (p^{\bigcup _{j<k-1}\L _j}_{H\cap\bigcup _{j<k-1}\L _j}\Rsh ^{L_{n_{k-1}}}_{H\cap L_{n_{k-1}}})$,
         if $k\geq 1$ and $H\subset \bigcup _{j<k}\L _j$.

\noindent
First we prove that the recursion works. By the duals of (sst1) and (sst2) $\CS$ has the unique top and $\Lev_0 (\CS )=\CS \cap L_{n_0}=\{ p \}$,
for some $p$. So for $p^\emptyset _\emptyset :=p$ and $\L _0 := \{p^\emptyset _\emptyset\}$ the sequence $\la \L _0\ra$ satisfies
($\L 1$) - ($\L 3$).

Suppose that a sequence $\la \L _0 , \dots ,\L _k\ra$ satisfies ($\L 1$) - ($\L 3$). Let $H\subset \bigcup _{j<k+1}\L _j$. Then
$H=(H\cap \bigcup _{j<k}\L _j) \cup (H\cap \L _k)$ and, by the assumption,
\begin{equation}\label{EQ2614}\textstyle
p^{\bigcup _{j<k}\L _j}_{H\cap \bigcup _{j<k}\L _j} \in \L _k \subset \CS \cap L_{n_k} \mbox{ and } H\cap \L _k \subset L_{n_k},
\end{equation}
\begin{equation}\label{EQ2616}\textstyle
H\cap \bigcup _{j<k}\L _j \subset \bigcup _{i<n_k}L_i
\end{equation}
and, by (L3) for $\CL _L$ we have  $p^{\bigcup _{j<k}\L _j}_{H\cap \bigcup _{j<k}\L _j}= q^{\bigcup _{i<n_k}L_i}_K$, for some
$K\subset \bigcup _{i<n_k}L_i$. So, by (\ref{EQ2614}) and Theorem \ref{T2627}(d) we have
$$
p^{\bigcup _{j<k}\L _j}_{H\cap \bigcup _{j<k}\L _j}\Rsh ^{L_{n_k}}_{H\cap L_{n_k}}
= q^{\bigcup _{i<n_k}L_i}_K\Rsh ^{L_{n_k}}_{H\cap L_{n_k}}
\in \Ip _{\la L, \leq _L \ra}(p^{\bigcup _{j<k}\L _j}_{H\cap \bigcup _{j<k}\L _j})\subset L_{n_k +1}
$$
and, by (sst3) the element
$p^{\bigcup _{j<k+1}\L _j}_H :=\pi _\CS (p^{\bigcup _{j<k}\L _j}_{H\cap\bigcup _{j<k}\L _j}\Rsh ^{L_{n_k}}_{H\cap L_{n_k}})$
is well defined and
belongs to $\CS \cap L_{n_{k+1}}$. Thus defining
$\L _{k+1} :=\{ p^{\bigcup _{j<k+1}\L _j}_H : H\subset \bigcup _{j<k+1}\L _j\}$
we have $\L _{k+1} \subset \CS \cap L_{n_{k+1}}$ and the sequence $\la \L _0 , \dots ,\L _k, \L _{k+1}\ra$ satisfies conditions
($\L 1$) - ($\L 3$). The recursion works indeed.

In order to prove that $\L$ is a copy of $R$ and $\CL _\L$ its labeling, using induction we show that
\begin{equation}\label{EQ2610'}\textstyle
\forall k\in \o \;\; \forall H\subset \bigcup _{j<k}\L _j \;\;p^{\bigcup _{j<k}\L _j}_H \in R^{\bigcup _{j<k}\L _j}_H .
\end{equation}
For $k=0$ we have $p^\emptyset _\emptyset \in R^\emptyset _\emptyset =R$.

Suppose that $k\in \o$ and
$p^{\bigcup _{j<k}\L _j}_H \in R^{\bigcup _{j<k}\L _j}_H$, for all $H\subset \bigcup _{j<k}\L _j$, and let
$H' \subset \bigcup _{j<k+1}\L _j$. Then $H'=(H' \cap \bigcup _{j<k}\L _j) \cup (H' \cap \L _k)$ and
we show that
\begin{equation}\label{EQ2642}\textstyle
p^{\bigcup _{j<k+1}\L _j}_{H'}\in R^{\bigcup _{j<k+1}\L _j}_{H'}= R^{\bigcup _{j<k}\L _j}_{H' \cap \bigcup _{j<k}\L _j} \cap R^{\L _k}_{H' \cap \L _k}.
\end{equation}
By ($\L 1$), ($\L 1$), and the induction hypothesis we have
\begin{equation}\label{EQ2618}\textstyle
L_{n_k}\ni p^{\bigcup _{j<k}\L _j}_{H' \cap \bigcup _{j<k}\L _j} \in R^{\bigcup _{j<k}\L _j}_{H' \cap \bigcup _{j<k}\L _j}
\end{equation}
thus, by (L3) and (L4) for $\CL _L$,
\begin{equation}\label{EQ2623}\textstyle
p^{\bigcup _{j<k}\L _j}_{H' \cap \bigcup _{j<k}\L _j}
=q^{\bigcup _{i<n_k}L_i}_K \in R^{\bigcup _{i<n_k}L_i}_K, \mbox{ where } K\subset \bigcup _{i<n_k}L_i .
\end{equation}
By (\ref{EQ2618}), (\ref{EQ2623}) and Lemma \ref{T2620}(a) we have $H' \cap \bigcup _{j<k}\L _j\cap \bigcup _{i<n_k}L_i =K \cap \bigcup _{j<k}\L _j$ and,
by ($\L 1$), $\bigcup _{j<k}\L _j\subset \bigcup _{j<k}L_{n_j}\subset \bigcup _{i<n_k}L_i$. Thus we have
\begin{equation}\label{EQ2634}\textstyle
H' \cap \bigcup _{j<k}\L _j =K \cap \bigcup _{j<k}\L _j .
\end{equation}
By ($\L 3$) and (\ref{EQ2623}) $p^{\bigcup _{j<k+1}\L _j}_{H'}
= \pi _\CS (p^{\bigcup _{j<k}\L _j}_{H' \cap\bigcup _{j<k}\L _j}\Rsh ^{L_{n_k}}_{H'\cap L_{n_k}})
= \pi _\CS (q^{\bigcup _{i\leq n_k}L_i}_{K \cup (H'\cap L_{n_k})})
=q^{\bigcup _{i< n_{k+1}}L_i}_F$, where $F\subset \bigcup _{i< n_{k+1}}L_i$ and, hence
$q^{\bigcup _{i< n_{k+1}}L_i}_F \leq _L q^{\bigcup _{i\leq n_k}L_i}_{K \cup (H'\cap L_{n_k})}$, which, by (\ref{EQ2605}), implies
\begin{equation}\label{EQ2641}\textstyle
F\cap (\bigcup _{i< n_k}L_i \cup L_{n_k}) = K \cup (H'\cap L_{n_k}) .
\end{equation}
Since $\bigcup _{j<k}\L _j\subset \bigcup _{i<n_k}L_i$ by  (\ref{EQ2634}) and (\ref{EQ2641}) we obtain
$F\cap \bigcup _{j<k}\L _j =K\cap \bigcup _{j<k}\L _j=H' \cap \bigcup _{j<k}\L _j $. Thus
\begin{equation}\label{EQ2643}\textstyle
p^{\bigcup _{j<k+1}\L _j}_{H'} =q^{\bigcup _{i< n_{k+1}}L_i}_F \in R ^{\bigcup _{i< n_{k+1}}L_i}_F \subset R^{\bigcup _{j<k}\L _j}_{F\cap \bigcup _{j<k}\L _j}
=R^{\bigcup _{j<k}\L _j}_{H' \cap \bigcup _{j<k}\L _j }.
\end{equation}
By ($\L 1$), $\L _k \subset L_{n_k}$, by (\ref{EQ2623}) we have $K\cap \L _k=\emptyset$ so, by (\ref{EQ2641}), $F\cap \L _k= H' \cap \L _k$ and, hence,
$
p^{\bigcup _{j<k+1}\L _j}_{H'}
=q^{\bigcup _{i< n_{k+1}}L_i}_F \in R ^{\bigcup _{i< n_{k+1}}L_i}_F \subset R^{\L _k}_{F\cap \L _k}=R^{\L _k}_{H'\cap \L _k},
$
which, together with (\ref{EQ2643})  gives  (\ref{EQ2642}).
So $\L \in \P (R)$ and $\CL _\L$ is a labeling of $\L$.

By Theorem \ref{T2627}, the labelings $\CL _L$ and $\CL _{\L}$ determine the reversed tree orderings $\leq _L$ and $\leq _\L$ on $L$ and $\L$
respectively, in the following way:
\begin{eqnarray}
q ^{\bigcup _{i<m}L_i}_{K '} \leq _L  q ^{\bigcup _{i<n}L_i}_{K''} & \Leftrightarrow &\textstyle m\geq n \;\land \;\; K' \cap \bigcup _{i<n}L_i =K'' ,
\label{EQ2606}\\
p ^{\bigcup _{j<k}\L_j}_{H '} \leq _\L  p ^{\bigcup _{j<l}\L_j}_{H''} & \Leftrightarrow &\textstyle \; k\geq l \;\;\land \;\; H' \cap \bigcup _{j<l}\L_j =H'' . \label{EQ2607}
\end{eqnarray}
and the sets $L_i$, $i\in \o$, and $\L _j$, $j\in \o$, are the corresponding levels.
In order to show that $\leq _\L \;=\; \leq _L \cap \;\L ^2$,
using induction we prove that for each $k\in \N$ we have
\begin{equation}\label{EQ2645}\textstyle
\forall u,v \in \bigcup _{j<k}\L _j \;\; (u <_\L v \Leftrightarrow u <_L v).
\end{equation}
For $k=1$ this follows from $|\L _0|=1$. Let $k\in \N$ and suppose that (\ref{EQ2645}) is true. We show that
\begin{equation}\label{EQ2646}\textstyle
\forall u,v \in \bigcup _{j<k}\L _j \cup \L _k \;\; (u <_\L v \Leftrightarrow u <_L v).
\end{equation}
Let $u,v \in \bigcup _{j<k}\L _j \cup \L _k$.

($\Rightarrow$) Let $u <_\L v$. If $u,v \in \bigcup _{j<k}\L _j$, then, by (\ref{EQ2645}), $u <_L v$. Otherwise, since
$\L _i$'s are the levels of the reversed tree $\la \L , \leq _\L \ra$, we have $u\in \L _k$ and $v\in \L _l$, for some $l<k$.
Also, there is $w\in \L _{k-1}$ such that $u<_\L w \leq _\L v$ and, since (\ref{EQ2645}) gives $w \leq _L v$ it remains to be shown that
$u<_L w$. By ($\L 3$)
$u= p^{\bigcup _{j<k}\L _j}_H =\pi _\CS (p^{\bigcup _{j<k-1}\L _j}_{H\cap\bigcup _{j<k-1}\L _j}\Rsh ^{L_{n_{k-1}}}_{H\cap L_{n_{k-1}}})
\leq _L p^{\bigcup _{j<k-1}\L _j}_{H\cap\bigcup _{j<k-1}\L _j}\Rsh ^{L_{n_{k-1}}}_{H\cap L_{n_{k-1}}}
<_L p^{\bigcup _{j<k-1}\L _j}_{H\cap\bigcup _{j<k-1}\L _j}$.
But, by  (\ref{EQ2607}) we have $u=p^{\bigcup _{j<k}\L _j}_H <_\L p^{\bigcup _{j<k-1}\L _j}_{H\cap\bigcup _{j<k-1}\L _j}\in \L _{k-1}$,
which implies that $p^{\bigcup _{j<k-1}\L _j}_{H\cap\bigcup _{j<k-1}\L _j}=w$ and, hence, $u <_L w$.

($\Leftarrow$) Let $u <_L v$. If $u,v \in \bigcup _{j<k}\L _j$, then, by (\ref{EQ2645}), $u <_\L v$. Otherwise,
since $\L _k \subset L_{n_k}$ and $L _n$'s are the levels of the reversed tree $\la L , \leq _L \ra$, we have $u\in \L _k$ and $v\in \L _l$, for some $l<k$.
Then, since $\CL _L$ and $\CL _\L$ are labelings of $L$ and $\L$,
\begin{equation}\label{EQ2647}\textstyle
u=p^{\bigcup _{j<k}\L _j}_{H'}=q^{\bigcup _{i< n_k}L_i}_F, \mbox{ where } H'\subset\bigcup _{j<k}\L _j \mbox{ and } F\subset \bigcup _{i< n_k}L_i ,
\end{equation}
\begin{equation}\label{EQ2648}\textstyle
v=p^{\bigcup _{j<l}\L _j}_{H''}=q^{\bigcup _{i< n_l}L_i}_G, \mbox{ where } H''\subset\bigcup _{j<l}\L _j \mbox{ and } G\subset \bigcup _{i< n_l}L_i ,
\end{equation}
so, by (\ref{EQ2607}), for a proof that $u\leq _\L v$ it remains to be shown that
$H' \cap \bigcup _{j< l}\L _j =H''$.

By (\ref{EQ2647}), (\ref{EQ2648}), Lemma \ref{T2620}(a) and (\ref{EQ2606}) we have
\begin{equation}\label{EQ2649}\textstyle
H' \cap \bigcup _{i< n_k}L_i =F \cap \bigcup _{j<k}\L _j\;\; \mbox{ and }\;\; H'' \cap \bigcup _{i< n_l}L_i =G \cap \bigcup _{j<l}\L _j,
\end{equation}
\begin{equation}\label{EQ2651}\textstyle
F\cap \bigcup _{i< n_l}L_i =G.
\end{equation}
Since $\bigcup _{j<l}\L _j \subset \bigcup _{i< n_l}L_i\subset \bigcup _{i< n_k}L_i$ from (\ref{EQ2649}) and (\ref{EQ2651}) we obtain
\begin{equation}\label{EQ2653}\textstyle
H' \cap \bigcup _{j<l}\L _j =F \cap \bigcup _{j<l}\L _j \;\; \mbox{ and }\;\; H'' \cap \bigcup _{j<l}\L _j =G \cap \bigcup _{j<l}\L _j,
\end{equation}
\begin{equation}\label{EQ2655}\textstyle
F\cap \bigcup _{j<l}\L _j =G\cap \bigcup _{j<l}\L _j.
\end{equation}
Now, since $H'' \cap \bigcup _{j<l}\L _j=H''$, the equality $H' \cap \bigcup _{j< l}\L _j =H''$ follows from (\ref{EQ2653}) and (\ref{EQ2655}).

The first equality in (\ref{EQ2613}) follows Theorem \ref{T2627}(c) applied
to $\L$, while the second follows from the equality $\leq _\L \;=\; \leq _L \cap \;\L ^2$.
\hfill $\Box$
\section{No splitting reals are added}\label{S10}
In this section we show that the poset $\P (R)$ shares one more property with the Sacks forcing. We recall that
if $\P$ is a forcing notion and $V_\P [G]$ a generic extension of the ground model $V$ by $\P$, then a real
$x\subset \o$ belonging to $V_\P [G]$ is called a {\it splitting real} iff $|A\cap x|=|A\setminus x|=\o$ for each infinite set
$A\subset \o$ belonging to $V$. It is well known that the Sacks forcing does not produce splitting reals and that the same holds for
the Miller rational perfect forcing (which does not have the Sacks property). Here we show that the poset $\P (R)$
(and, consequently, the first iterand $\P$ in the two-step iteration $\P \ast \pi$, see Section \ref{S6}) has this property as well.
\begin{te}  \label{T2626}
The forcing $\la \P (R), \subset \ra$ does not produce splitting reals.
\end{te}
\dok
We prove that for each $\P (R)$-name $\t$
\begin{equation}\label{EQ2600}\textstyle
R \Vdash \t \subset \check{\o } \Rightarrow \exists S\in (( [\o ]^\o)^V)\check{\;}\;\; (S\subset \t \lor S\subset \check{\o }\setminus \t).
\end{equation}
Thus, working in $V$ and assuming that $\P (R) \ni A\Vdash \t  \subset \check{\o }$ it is sufficient to find $\Lambda \in \P (A)$ and
$S\in [\o ]^\o$ such that
\begin{equation}\label{EQ2601}\textstyle
\Lambda  \Vdash   \check{S}\subset \t \;\;\lor \;\; \Lambda  \Vdash \check{S}\subset \check{\o }\setminus \t .
\end{equation}
Since the sets $\D _n=\{ D\in \P (R): D\Vdash \check{n}\in \t \; \lor \; D\Vdash\check{n}\not\in \t \}$, $n\in \o$,
are dense in $\P (R)$, by Theorem \ref{T2637} the set of fusions of the sequence $\la D_n :n\in \o \ra$ is dense as well and, hence,
there is a fusion $L=\bigcup _{n\in \o}L_n \subset A$.
So we have  $L_n =\{ q ^{\bigcup _{i<n}L_i}_K : K\subset \bigcup _{i<n}L_i \}$, where $q ^{\bigcup _{i<n}L_i}_K \in  L^{\bigcup _{i<n}L_i}_K$,
and, by (\ref{EQ2602}), for each $n\in \o$ and each $K \subset \bigcup _{i<n}L_i$ there is $D\in \D _n$ such that
$L^{\bigcup _{i<n}L_i}_K \subset D$. Thus
\begin{equation}\label{EQ2604}\textstyle
\forall n\in \o \;\; \forall K\subset \bigcup _{i<n}L_i\;\;
(L^{\bigcup _{i<n}L_i}_K\Vdash \check{n}\in \t \;\; \lor \;\;  L^{\bigcup _{i<n}L_i}_K\Vdash \check{n} \not\in \t ).
\end{equation}
By Theorem \ref{T2627} $\la L, \leq _L \ra$ is a reversed tree and for each $K\subset \bigcup _{i<n}L_i$ we have
$L^{\bigcup _{i<n}L_i}_K= (-\infty , q^{\bigcup _{i<n}L_i}_K]_{\la L, \leq _L \ra}$. So, by (\ref{EQ2604}), $L=L' \cup L''$ is a coloring of $L$, where
$$
L'=\{ q^{\bigcup _{i<n}L_i}_K \in L : (-\infty , q^{\bigcup _{i<n}L_i}_K]_{\la L, \leq _L \ra} \Vdash \check{n}\in \t \} ,
$$
$$
L''=\{ q^{\bigcup _{i<n}L_i}_K \in L : (-\infty , q^{\bigcup _{i<n}L_i}_K]_{\la L, \leq _L \ra} \Vdash \check{n}\not\in \t \} .
$$
Now by Theorem \ref{T2636}
there is a monochromatic strong reversed subtree $\CS$ of the reversed tree $\la L , \leq _L\ra$. Let $S=\{ n_k :k\in \o \}$ be the level
set of $\CS$.

First suppose that $\CS \subset L'$.
By Theorem \ref{T2629} there is a copy $\Lambda =\bigcup _{k\in \o}\L _k \subset \CS $ such that
$\L _k=\{ p ^{\bigcup _{j<k}\L _j}_H : H\subset  \bigcup _{j<k}\L _j\} \subset L_{n_k}$ and $p ^{\bigcup _{j<k}\L _j}_H \in \L ^{\bigcup _{j<k}\L _j}_H$.

We prove that $\Lambda \Vdash \check{S}\subset \t$, that is
$\Lambda \Vdash \check{n_k}\in \t $, for all $k\in \o$.
Since $\{ \L ^{\bigcup _{j<k}\L _j}_H :H\subset \bigcup _{j<k}\L _j\}$ is an antichain in $\la \P (R),\subset \ra$ maximal below $\L$, for a proof of
$\Lambda \Vdash \check{n_k}\in \t $
it is sufficient to show that for each $H\subset  \bigcup _{j<k}\L _j$ we have
\begin{equation}\label{EQ2656}\textstyle
 \L ^{\bigcup _{j<k}\L _j}_H\Vdash \check{n_k}\in \t .
\end{equation}
By Theorem \ref{T2629} we have
$\L ^{\bigcup _{j<k}\L _j}_H =\L \cap (-\infty , p ^{\bigcup _{j<k}\L _j}_H ]_{\la L , \leq _L  \ra} $ so
\begin{equation}\label{EQ2657}\textstyle
\L ^{\bigcup _{j<k}\L _j}_H
\subset (-\infty , p ^{\bigcup _{j<k}\L _j}_H ]_{\la L , \leq _L  \ra} .
\end{equation}
Since $p ^{\bigcup _{j<k}\L _j}_H\in L_{n_k}$ we have
$p ^{\bigcup _{j<k}\L _j}_H= q ^{\bigcup _{i<n_k}L_i}_K$, for some $K\subset \bigcup _{i<n_k}L_i$.
Now $\L \subset \CS \subset L'$ implies $q ^{\bigcup _{i<n_k}L_i}_K\in L'$; thus
$(-\infty , q^{\bigcup _{i<n_k}L_i}_K]_{\la L , \leq _L  \ra} \Vdash \check{n_k}\in \t$
and, by (\ref{EQ2657}), $\L ^{\bigcup _{j<k}\L _j}_H \Vdash \check{n_k}\in \t$.
So (\ref{EQ2656}) is proved.

If $\CS \subset L''$, then in a similar way we prove that $\Lambda \Vdash \check{S}\subset \check{\o }\setminus \t$.
\kdok

\noindent
{\bf Acknowledgement.}
Research was supported by the Ministry of Education and Science of the Republic of Serbia (Project 174006) and grants from CNRS and NSERC.

\end{document}